\documentclass[twoside]{article}
\usepackage{mathrsfs}
\usepackage{graphicx,amssymb,mathrsfs,amsmath}
\catcode`@=11
\long\def\@makefntext#1{\noindent #1}
\newskip\tabcentering \tabcentering=1000pt plus 1000pt minus 1000pt
\def\REF#1{\par\hangindent\parindent\indent\llap{#1\enspace}\ignorespaces} 
\def\MCH#1#2{\setbox0=\hbox{\raise#1\hbox{#2}}\smash{\box0}}

\def\@evenfoot{}\def\@oddfoot{}

\def\@oddhead{\hbox to \textwidth{\footnotesize{\it
Asymptotic behavior for density estimators} \hfill\thepage}}

\def\sec#1{\vskip 3mm\leftline{\bf #1}\vskip 1mm}

\catcode`@=12

\def\bc{\begin{center}}
\def\ec{\end{center}}
\def\no{\noindent}
\def\hang{\hangindent\parindent}
\def\textindent#1{\indent\llap{\qquad #1\ \ \enspace}\ignorespaces}
\def\ref{\par\hang\textindent}

\begin{document}
\abovedisplayskip=6pt plus 1pt minus 1pt \belowdisplayskip=6pt
plus 1pt minus 1pt
\thispagestyle{empty} \vspace*{-1.0truecm} \noindent
\vskip 10mm

\bc{\large\bf  Limit theorems for kernel density estimators  under dependent samples
 \footnotetext{\footnotesize \\
Supported by the National Natural Science Foundation of China (Grant
Nos. 11171303, 61273093) and the Specialized Research Fund for the
Doctor Program
  of Higher Education (Grant No. 20090101110020).}} \ec

\vskip 2mm
\bc{\bf  Yuexu Zhao, Zhengyan Lin\\
{\small\it  Department  of Mathematics, Zhejiang University,
Hangzhou {\rm 310027},  China}}\ec

\vskip 1 mm

\noindent{\small {\small\bf Abstract.} \textmd{
 In this paper, we construct a moment inequality for mixing dependent random variables, it is of independent interest. As applications,
  the consistency of the kernel density estimation  is investigated.
 Several  limit theorems are established: First, the central limit theorems for
 the kernel density estimator $f_{n,K}(x)$  and its distribution function are constructed.
 Also, the  convergence rates of
 $\|f_{n,K}(x)-Ef_{n,K}(x)\|_{p}$ in sup-norm loss
  and integral $L^{p}$-norm loss are proved.
 Moreover, the a.s. convergence rates of the supremum of $|f_{n,K}(x)-Ef_{n,K}(x)|$ over a compact set and the whole real line
 are obtained.
 It is showed, under suitable  conditions on the mixing rates, the kernel
function and the bandwidths,  that the optimal rates for i.i.d.
random variables are also optimal for dependent ones.}
 \vspace{1mm}\baselineskip 12pt

\no{\small\bf Keywords:} Kernel density estimator;
consistency; convergence rate;  mixing rate\ \

\no{\small\bf Mathematics  Subject Classification (2000):}\ \ {\rm
Primary 62G07, 60G10; Secondary  60F15.}}

\sec{1.\quad Introduction}
Let $X, X_{1},X_{2},...$ be independent and identically distributed
(i.i.d.) random variables with common density $f$, and $K$ be a
bounded integrable kernel (a measurable function on $\mathbb{R}$), the classical kernel
density estimators (KDEs) of $f$ based on the observations
$X_{1},...,X_{n}$ are defined as
$$f_{n,K}(x)=\frac{1}{nh_{n}}\sum_{i=1}^{n}K\Big{(}\frac{X_{i}-x}{h_{n}}\Big{)},~~~~x\in \mathbb{R},\eqno(1.1)$$
where the bandwidths $\{h_{n}, n\geq1\}$ satisfy some regularity
conditions.

Since the famous work done by Rosenblatt [25] and Parzen [18],
the limit behavior for the  KDEs  has become an active subject.
 For the case of i.i.d. data,
 see, for example,
 Bickel and Rosenblatt [1], Silverman [34] and
Stute [35, 36]. Using empirical process approach, Einmahl and Mason
[4, 5] studied the uniform consistency and uniform consistency in
bandwidth, respectively. Gin\'{e} and Guillou [8, 9]
 investigated  the exact rates of almost sure (a.s.)
convergence  of the supremum  over adaptive intervals and all of
$\mathbb{R}^{d}$, and Gin\'{e}, Koltchinskii and Zinn [10] obtained
weighted uniform consistency of  KDEs, and so forth. As to weakly
dependent observations, F\"{o}ldes [6], R\"{u}schendorf [27], Sarda
and Vieu [28], Peligrad [20] and Liebscher [13] studied the strong
convergence of density estimators for $\phi-$mixing samples.
Rosenblatt [25],
 Nze and Rios [17], Liebscher [15] investigated  a.s.
convergence of kernel estimators for $\alpha-$mixing random
variables. For other results, one can  refer to Neumann [16], Woodroofe [37, 38],
Wu et al. [39], Yakowitz [40], and the reference therein.
 However, most of the work mentioned as above on  a.s. convergence rates in sup-norm loss
  under dependent data are not optimal.  Yu [41] obtained the  best possible minimax rates for stationary sequences satisfying certain $\beta-$mixing
conditions at the cost of sufficient  smoothness  for density functions.
  The purpose of the present article is to investigate the   consistency of
  the KDEs, and tries to get the optimal convergence rates for certain dependent observations.  More precisely, we require the random variables to be $\rho-$mixing,
  which  is defined as follows:
  \\
  \\
  {\bf Definition.}\textit{
  Suppose that $X_{1},X_{2},...$ is a sequence of random variables on a probability space
   $(\Omega,\mathscr{F},P)$. Set  $\mathscr{F}_{n}^{-}=\sigma(X_{i}, 1\leq i\leq n)$,
   $\mathscr{F}_{n}^{+}=\sigma(X_{i}, i\geq n),$ define
   $$\rho(n)=\sup_{k\geq1}\sup_{X\in L^{2}(\mathscr{F}_{k}^{-})}\sup_{Y\in L^{2}(\mathscr{F}_{k+n}^{+})}
              \frac{|EXY-EXEY|}{\sqrt{E(X-EX)^{2}E(Y-EY)^{2}}}.\eqno(1.2)$$
  The sequence  $X_{1},X_{2},...$ is said to be $\rho-$mixing if $\rho(n)\rightarrow0$ as
  $n\rightarrow\infty$.}

   This definition  was introduced by Kolmogorov and Rozanov [12].
  As is known,  the asymptotic behavior of $\rho-$mixing sequences have  received much well-deserved attention,
  and a variety of  elegant  results  have been obtained.  See, for example, Lin and Lu [15],  Peligrad [19-21], Peligrad and Shao [23], Peligrad and  Utev  [24], Shao [29-33], and so forth.

  Let $X_{1},X_{2},...$ be a sequence of  stationary $\rho-$mixing
  random variables with  density $f$. replace the independent observations by the
 $\rho-$mixing ones in (1.1), one gets the corresponding density estimator of $f$ for the dependent random variables.

  In this  article, we  devote ourselves  to doing  three things.
  The first one is to study convergence in distribution of the  estimator $f_{n,K}(x)$ both as an estimation for the true density function $f(x)$ and as an estimation
   $F_{n,K}(x)=\int_{-\infty}^{x}f_{n,K}(t)dt$ for the true distribution function $F(x)$ of $X$. The second  is to
   investigate the convergence rates  for the difference  of $f_{n,K}(x)$ and its mean in  sup-norm loss
  and integral $L^{p}$-norm loss. Our third goal is to discuss
  the strong uniform   convergence rates of $|f_{n,K}(x)-Ef_{n,K}(x)|$  over a compact set of $\mathbb{R}$ and  the whole real line $\mathbb{R}$, respectively.
  Of course,
 a natural question is posed as follows:
  Whether the optimal convergence rates could be achieved?
 The answer is affirmative
  for  i.i.d. observations. As is known
 a variety of sharp results  have been established, see, for example,
  Einmahl and Mason [4, 5],  Gin\'{e} and Guillou [8, 9],   Gin\'{e}, Koltchinskii and Zinn
  [10]. However, that in general is not the case for
  dependent samples. To obtain the best possible convergence rates,  some different methods from those for i.i.d. case should be developed.  The
  present paper tries to do this. Our technical proofs consist in applications of the blocking
techniques,
  the martingale methods  and some inequalities.
  It is showed that the optimal convergence rates for i.i.d. random
variables are also optimal for dependent ones.

 The remainder of the paper is structured as follows. Section 2 introduces some notation and assumptions.
 Section 3 formulates several  results  on the weak convergence. Section 4 constructs the rates of  $\|f_{n,K}(x)-Ef_{n,K}(x)\|_{p}$ in
  the sup-norm loss and integral $L^{p}$-norm loss, while Section 5 derives  the rates of strong  uniform consistency for  KDEs.
 Some useful results are stated in the Appendix.

\sec{2.\quad Notation and assumptions}
   In this section, we present  some basic notation and assumptions which will be used in the sequel.
    Let $X, X_{1}, X_{2},...$ be a sequence of non-degenerated and stationary $\rho-$mixing random variables. Denote $K_{i}(x)=K\big{(}(X_{i}-x)/h_{n}\big{)}$ for fixed $n\in \mathbb{N}$,  where $K$  is a
   measurable function satisfying some regularity conditions.
    $f(x)$ is the unknown  density function of $X$  with respect to Lebesgue measure.
    For Borel measurable function $g$ and Borel measure $\mu$, let $L^{p}:=L^{p}(\mu)$ be the usual Lebesgue spaces of real-valued functions  normed  by  $\|\cdot\|_{p}$.
    As usual, write $\|g\|_{p}=\big{(}\int_{\mathbb{R}}|g(x)|^{p}d\mu(x)\big{)}^{1/p}$ for $1\leq
   p<\infty.$
   Define for any nonnegative
integer $s$ the spaces $\mathcal{C}^{s}(\mathbb{R})$ of all bounded continuous real-valued functions
that are $s-$times continuously differentiable on $\mathbb{R}$.
 $\textmd{I}(\cdot)$ is the indicator function.  $[z]$ denotes the integer part of $z$, $\log x=\log (x\vee e)$.
 $a_{n}=O( b_{n})$ means $\limsup_{n\rightarrow\infty}a_{n}/b_{n}<\infty,$
  $a_{n}=o( b_{n})$ stands for $\limsup_{n\rightarrow\infty}a_{n}/b_{n}=0$,
  and $a_{n}\asymp b_{n}$ means $0<\liminf_{n\rightarrow\infty}a_{n}/b_{n}\leq\limsup_{n\rightarrow\infty}a_{n}/b_{n}<\infty$.  The letter $C$
   with subscripts denotes some finite and positive universal
   constants, it may take different values in each appearance.
   \\
    \\
    Some assumptions are formulated  below:
    \\
   (B1) $h_{n}\searrow0$ and $ nh_{n}\rightarrow\infty$ as $n\rightarrow\infty$.
   \\
    (B2) $h_{n}\asymp n^{-\delta}l(n)$ for  $0<\delta\leq1$,
         $l(n)$ is a slowly varying function.
    \\
   (C1) the density function $f(x)$  of $X$ is  uniformly bounded on $\mathbb{R}$.
   \\
    (C2) the density function $f(x)$  of $X$ is uniformly continuous and uniformly bounded on $\mathbb{R}$.
    \\
   (K1) $K$ is a real-valued measurable function satisfying  $\sup_{x\in\mathbb{R}}|K(x)|<\infty$ and  $\int_{-\infty}^{\infty}|K(x)|dx<\infty$.
   \\
   (K2) $K$ is  a real-valued measurable function with  compact support on  $\mathbb{R}$, and  satisfies  Lipschitz condition.
   \\
   \\
  {\bf Remark 1.} Condition (B2) is a little more stronger than  (B1). in other words, (B2) does not allow the bandwidths $h_{n}$ to go to zero very slowly as $n\rightarrow\infty$. For example, the form of the bandwidths such as $h_{n}=1/(\log n)^{p}$, for all $p>0$, is excluded. But we would like to point out that most of the  bandwidths including the optimal ones are contained in (B2).

\sec{3. \quad Central limit theorems for  KDEs and their distribution functions}
  Consider the KDE $f_{n, K}(x)$ defined in (1.1). The  aim  of  this section is to
 investigate the CLT  for
 $f_{n, K}(x)$ and $F_{n, K}(x)$. The classical theory of this subject was developed mostly in the
 1950s, and it is an important theory in probability and statistics.
 Our first result reads  as follows:
\\
\\
{\bf Theorem 3.1.} \textit{Suppose that conditions (B1), (C2) and (K1)
hold.  Further assume that $f(x)>0$ and  $\sum_{i=0}^{\infty}\rho(2^{i})<\infty$.
 Then we have
$$\sqrt{nh_{n}}(f_{n, K}(x)-E f_{n, K}(x))\stackrel{d}{\longrightarrow} N(0, \|K\|_{2}^{2}f(x)),\eqno(3.1)$$
where $"\stackrel{d}{\longrightarrow}"$ stands for convergence in
distribution.}
 \\
 \textit{Proof.}
 For any fixed $x\in\mathbb{R}$, we  use the following decomposition:
 $$f_{n, K}(x)-E f_{n, K}(x)=\big{[}f_{n, K}(x)-\Gamma_{n,K}(x)\big{]}
 +\big{[}\Gamma_{n,K}(x)-Ef_{n, K}(x)\big{]},\eqno(3.2)$$
where
$$\Gamma_{n,K}(x):=\frac{1}{nh_{n}}\sum_{i=1}^{n}E\big{[}K((X_{i}-x)/h_{n})|\mathscr{F}_{i-1}\big{]},~~\mathscr{F}_{i}:=\sigma(X_{j}, j\leq i),~\mathscr{F}_{0}=\{\emptyset, \Omega\}.$$
Thus,  (3.1) will be derived if one can show  that
$$\Gamma_{n,K}(x)-Ef_{n, K}(x)=o_{P}\bigg{(}\frac{1}{\sqrt{nh_{n}}}\bigg{)}\eqno(3.3)$$
and
$$\sqrt{nh_{n}}\big{(}f_{n, K}(x)-\Gamma_{n,K}(x)\big{)}\stackrel{d}{\longrightarrow} N(0, \|K\|_{2}^{2}f(x)).\eqno(3.4)$$

 We first prove (3.3). However, some preliminary work is needed.  Denote $\mathbb{N}^{+}=\{1,2,...\}$, and let $\mathbb{I}_{k}$ be the integer interval $[2^{k}, 2^{k+1})$. Clearly, for each $n\in\mathbb{N}^{+}$, there exists
 integer $k_{n}\geq0$ such that $2^{k_{n}}\leq n<2^{k_{n}+1}$.
 Moreover, for  $0<\beta<\alpha<1$,  let
   $p_{k}=[2^{\alpha k}], q_{k}=[2^{\beta k}], r_{k}=[2^{k}/(p_{k}+q_{k})]$,  then the integer set  $\mathbb{I}_{k}$ can be blocked as follows:\\
 $~~~~~~~~~~~~~~~~~~~~~~~\mathbb{I}_{k}(m)=[2^{k}+(m-1)(p_{k}+q_{k}), 2^{k}+(m-1)q_{k}+mp_{k})\cap \mathbb{N}^{+},$\\
 $~~~~~~~~~~~~~~~~~~~~~~~\mathbb{J}_{k}(m)=[ 2^{k}+(m-1)q_{k}+mp_{k}, 2^{k}+m(p_{k}+q_{k}))\cap \mathbb{N}^{+},$\\
 $~~~~~~~~~~~~~~~~~~~~~~~1\leq m\leq r_{k},~~\mathbb{J}_{k}(r_{k}+1)=[2^{k}+r_{k}(p_{k}+q_{k}),  2^{k+1})\cap \mathbb{N}^{+}.$\\
 It is easy to see that $r_{k}\sim 2^{(1-\alpha)k}$. According to the  symbols as above, there
 also
 exists some integer $m_{n}\geq0$ such that $n\in\mathbb{I}_{k_{n}}(m_{n})\cup\mathbb{J}_{k_{n}}(m_{n})$. For simplicity, we introduce some extra notation as follows: Denote
 $$W_{k}^{m}(x):=\sum_{i\in \mathbb{I}_{k}(m)}\Big{[}E(K_{i}(x)|\mathscr{F}_{i-1})-EK_{i}(x)\Big{]},
 ~~V_{k}^{m}(x):=\sum_{i\in \mathbb{J}_{k}(m)}\Big{[}E(K_{i}(x)|\mathscr{F}_{i-1})-EK_{i}(x)\Big{]}.$$
 Then it follows that
 \begin{eqnarray*}
 &&nh_{n}\big{[}\Gamma_{n,K}(x)-Ef_{n, K}(x)\big{]}=\Bigg{[}\sum_{k=0}^{k_{n}-1}\sum_{m=1}^{r_{k}}W_{k}^{m}(x)+\sum_{m=1}^{m_{n}-1}W_{k_{n}}^{m}(x)\Bigg{]}\\
  &&~~~+\Bigg{[}\sum_{k=0}^{k_{n}-1}\sum_{m=1}^{r_{k}+1}V_{k}^{m}(x)+\sum_{m=1}^{m_{n}-1}V_{k_{n}}^{m}(x)\Bigg{]}
+\sum_{i=N_{n}}^{n}\Big{[}E(K_{i}(x)|\mathscr{F}_{i-1})-E
K_{i}(x)\Big{]},
  \end{eqnarray*}
 where $N_{n}=2^{k_{n}}+(m_{n}-1)(p_{k_{n}}+q_{k_{n}})$.

 Thus, in order to verify (3.3), it suffices to show that the sums on the big
 blocks satisfy
  $$E\Bigg{[}\sum_{k=0}^{k_{n}-1}\sum_{m=1}^{r_{k}}W_{k}^{m}(x)+\sum_{m=1}^{m_{n}-1}W_{k_{n}}^{m}(x)\Bigg{]}^{2}=o\big{(}nh_{n}\big{)}.\eqno(3.5)$$
Note that the left-hand side of (3.5) is controlled by
$$2E\Bigg{[}\sum_{k=0}^{k_{n}-1}\sum_{m=1}^{r_{k}}W_{k}^{m}(x)\Bigg{]}^{2}
+2E\Bigg{[}\sum_{m=1}^{m_{n}-1}W_{k_{n}}^{m}(x)\Bigg{]}^{2}=:\Sigma_{1}+\Sigma_{2}.\eqno(3.6)$$
So we only need to show that
$$\Sigma_{1}=o\big{(}nh_{n}\big{)},~~~~\Sigma_{2}=o\big{(}nh_{n}\big{)}.\eqno(3.7)$$
Using the towering property and  Jensen$^{,}$s inequality for the conditional expectations together with Lemma A.4 with $p=2$, we have
 \begin{eqnarray*}
 &&\Sigma_{1}=2E\Bigg{[}\sum_{k=0}^{k_{n}-1}E\Bigg{(}\sum_{m=1}^{r_{k}}W_{k}^{m}(x)\Big{|}\mathscr{F}_{2^{k}+r_{k}-1}\Bigg{)}\Bigg{]}^{2}\\
 &&~~~~\leq C\log^{2}(2k_{n})\sum_{k=0}^{k_{n}-1}\rho^{2}(q_{k})E\Bigg{[}\sum_{m=1}^{r_{k}}W_{k}^{m}(x)\Bigg{]}^{2}\\
 &&~~~~=C\log^{2}(2k_{n})\sum_{k=0}^{k_{n}-1}\rho^{2}(q_{k})E\Bigg{[}\sum_{m=1}^{r_{k}}E\Big{(}W_{k}^{m}(x)\big{|}\mathscr{F}_{t_{k}(m)}\Big{)}\Bigg{]}^{2}\\
 &&~~~~\leq C\log^{2}(2k_{n})\sum_{k=0}^{k_{n}-1}\rho^{4}(q_{k})\log^{2}(2r_{k})
  \sum_{m=1}^{r_{k}}E\Big{(}W_{k}^{m}(x)\Big{)}^{2}\\
  &&~~~~\leq C\log^{2}(2k_{n})\sum_{k=0}^{k_{n}-1}2^{k}\rho^{4}(q_{k})\log^{2}(2r_{k})\log^{2}(2p_{k})\|K_{1}(x)\|_{2}^{2},
 \end{eqnarray*}
where
$t_{k}(m)=2^{k}+(-q_{k-1})\textrm{I}(m=1)+(m-2)(p_{k}+q_{k})\textrm{I}(m\neq1)+p_{k},
1\leq m\leq r_{k}$.

 Recalling the condition imposed on the mixing rates,  without loss of generality (w.l.o.g.), suppose that $\rho(n)\leq 1/\log n$, and observe that
$$EK_{1}^{2}(x)=h_{n}\int_{-\infty}^{\infty}K^{2}(u)f(x+h_{n}u)du\leq h_{n}\|f\|_{\infty}\|K\|_{\infty}\|K\|_{1}.\eqno(3.8)$$
Then, applying Lemma A.2, we can get
 $$\Sigma_{1}\leq  C\log^{2}(2k_{n})EK_{1}^{2}(x)\sum_{k=0}^{k_{n}-1}2^{k}k^{-2}
 \leq Cnh_{n}(\log\log n)^{2}(\log n)^{-1}=o(nh_{n}).\eqno(3.9)$$
Similarly, we have
$$\Sigma_{2}
\leq Cnh_{n}(\log\log n)^{2}(\log n)^{-2}=o(nh_{n}).\eqno(3.10)$$
Combining (3.9) and (3.10) yields (3.3).

 As to (3.4), note that
 $$\sqrt{nh_{n}}\big{(}f_{n, K}(x)-\Gamma_{n,K}(x)\big{)}=\frac{1}{\sqrt{nh_{n}}}\sum_{i=1}^{n}\Big{[}K_{i}(x)
 -E\big{(}K_{i}(x)|\mathscr{F}_{i-1}\big{)}\Big{]}.\eqno(3.11)$$
We next estimate  the conditional variance
\begin{eqnarray*}
&&\frac{1}{nh_{n}}\sum_{i=1}^{n}E\Big{\{}\big{[}K_{i}(x)
 -E\big{(}K_{i}(x)|\mathscr{F}_{i-1}\big{)}\big{]}^{2}|\mathscr{F}_{i-1}\Big{\}}\\
 &&~~=\frac{1}{nh_{n}}\sum_{i=1}^{n}E\big{[}K_{i}^{2}(x)|\mathscr{F}_{i-1}\big{]}
 -\frac{1}{nh_{n}}\sum_{i=1}^{n}\big{[}E\big{(}K_{i}(x)|\mathscr{F}_{i-1}\big{)}\big{]}^{2}=:\Xi_{1}-\Xi_{2}.
\end{eqnarray*}
For $\Xi_{1}$, observe that
$$P\big{(}|\Xi_{1}-\|K\|_{2}^{2}f(x)|>\epsilon\big{)}
=P\big{(}|\Xi_{1}-h_{n}^{-1}EK_{1}^{2}(x)+h_{n}^{-1}EK_{1}^{2}(x)-\|K\|_{2}^{2}f(x)|>\epsilon\big{)}.\eqno(3.12)$$
Clearly,  on account of   condition (C2) and  Bonchner$^{,}$s lemma,
we have
$$h_{n}^{-1}EK_{1}^{2}(x)
=\int_{-\infty}^{\infty}K^{2}(u)f(x+h_{n}u)du\rightarrow\|K\|_{2}^{2}f(x).\eqno(3.13)$$
Therefore, by Lemma A.2 and Jensen$^{,}$s inequality together with condition (B1), the right-hand side of (3.12), for large $n$, is controlled by
$$P\big{(}|\Xi_{1}-h_{n}^{-1}EK_{1}^{2}(x)|>\epsilon/2\big{)}
\leq C\epsilon^{-2}n^{-1}h_{n}^{-2}EK_{1}^{4}(x)\leq
C\epsilon^{-2}(nh_{n})^{-1}\|f\|_{\infty}\rightarrow0.\eqno(3.14)$$
For $\Xi_{2}$,  we have with probability one,
$$h_{n}^{-1}E\big{(}K_{i}(x)|\mathscr{F}_{i-1}\big{)}
=\int_{-\infty}^{\infty}K(u)f(x+h_{n}u|\mathscr{F}_{i-1})du\leq
\|K\|_{\infty}\|f\|_{\infty}.$$ Then, it follows that
$$E\Xi_{2}\leq h_{n}\|K\|_{\infty}\|f\|_{\infty}\rightarrow0.\eqno(3.15)$$
Combining (3.12)$-$(3.15) yields
$$ \frac{1}{nh_{n}}\sum_{i=1}^{n}E\Big{\{}\big{[}K_{i}(x)
 -E\big{(}K_{i}(x)|\mathscr{F}_{i-1}\big{)}\big{]}^{2}|\mathscr{F}_{i-1}\Big{\}}
 \stackrel{P}{\longrightarrow}\|K\|_{2}^{2}f(x).\eqno(3.16)$$
Moreover, applying the $C_{r}$-inequality, we have
 $$E\big{(}K_{i}(x)
 -E\big{(}K_{i}(x)|\mathscr{F}_{i-1}\big{)}\big{)}^{2}\leq 4EK_{i}^{2}(x)
 \leq 4h_{n}\|f\|_{\infty}.\eqno(3.17)$$
 Thus on account of (B1), the Lindeberg  condition
 \begin{eqnarray*}
 && \frac{1}{nh_{n}}\sum_{i=1}^{n}E\Big{[}K_{i}(x)
 -E\big{(}K_{i}(x)|\mathscr{F}_{i-1}\big{)}\Big{]}^{2}\textmd{I}(|K_{i}(x)
 -E(K_{i}(x)|\mathscr{F}_{i-1})|>\epsilon\sqrt{nh_{n}})\\
 &&~~\leq 4h_{n}^{-1}EK_{1}^{2}(x)\textmd{I}(|K_{1}(x)|>\epsilon\sqrt{nh_{n}}/2)
 \leq 4\|f\|_{\infty}\int_{|K(u)|>\epsilon\sqrt{nh_{n}}/2}K(u)du=o(1)
 ~~~~(3.18)
 \end{eqnarray*}
 holds for any $\epsilon>0$.

Finally,  according to (3.16) and (3.18), then using the martingale central limit theorem together with  Slutsky$^{,}$s  theorem gives (3.1).
\\
\\
{\bf Remark 2.}  Let us consider  the
deviation of the kernel density estimator with respect to the true density function.  Note that
$$f_{n,K}(x)-f(x)=\big{[}f_{n,K}(x)-Ef_{n,K}\big{]}+\big{[}Ef_{n,K}(x)-f(x)\big{]}.\eqno(3.19)$$
   The first term on the right-hand side of (3.19) is the probabilistic term,  while the second term is the bias.
   If  (C2) and the conditions imposed on the kernel $K$ in Theorem 3.1 are replaced by
   \\
   \\
 (C3) the density function $f(x)$ is uniformly  bounded, $f(x)\in \mathcal{C}(\mathbb{R})$
       and $\sup_{x} |f^{'}(x)|<\infty$,
  \\
 (K3) $K$  satisfies  $\sup_{x}|K(x)|<\infty, ~\int_{-\infty}^{\infty}|xK(x)|dx<\infty$,
\\
 then  applying Taylor$^{,}$s expansion, we have for some $0<\upsilon<1$,
 \begin{eqnarray*}
&&\big{|}Ef_{n,K}(x)-f(x)\big{|}=\bigg{|}\int_{-\infty}^{\infty}K(y)\big{[}f(x-h_{n}y)-f(x)\big{]}dy\bigg{|}\\
&&~~~~~~~~~~~~~~~~~~~~~~~~=h_{n}\bigg{|}\int_{-\infty}^{\infty}y
K(y)f^{'}(x-\upsilon h_{n}y)dy\bigg{|}
 =O(h_{n}).
 \end{eqnarray*}
 Therefore, (3.1) holds whenever $nh_{n}^{3}\rightarrow0$ as $n\rightarrow\infty$.

In fact,  the bias can always be balanced with the probabilistic term
  by calibrating the
normalizing sequence $\{h_{n}, n\geq1\}$, provided enough regularity for $K$ and $f$ are assumed.

  \vspace{3mm}
  Another interesting problem  is  the weak convergence for the distribution
function of the KDE. More precisely, denote
$F_{n,K}(x)=\int_{-\infty}^{x}f_{n,K}(t)dt$, we construct the CLT
for the difference between $F_{n,K}(x)$ and its mean.
\\
  \\
  {\bf Theorem 3.2.} \textit{Suppose that condition (B1) holds, and that $\sum_{i=0}^{\infty}\rho(2^{i})<\infty$. Further, assume that the density function $f(x)$ is continuous
   and positive on $\mathbb{R}$, and that the kernel $K$ is symmetric and $\int_{-\infty}^{\infty}K(x)dx=1$.
 Then we have
 $$\sqrt{n}(F_{n,K}(x)-EF_{n,K}(x))\stackrel{d}{\longrightarrow} N(0, F(x)(1-F(x))),\eqno(3.20)$$
 where $F(x)$ is the true distribution function of $X$.}
 \\
 \textit{Proof.} Observe that
 $$F_{n,K}(x)-EF_{n,K}(x)=[F_{n,K}(x)-\Lambda_{n,K}(x)]+[\Lambda_{n,K}(x)-EF_{n,K}(x)],\eqno(3.21)$$
 where
 $$\Lambda_{n,K}(x):=\frac{1}{nh_{n}}\sum_{i=1}^{n}E\bigg{[}\int_{-\infty}^{x}K_{i}(t)dt\Big{|}\mathscr{F}_{i-1}\bigg{]}.$$

 Along the similar proof lines as those  of (3.3), one can get
 $$\Lambda_{n,K}(x)-EF_{n,K}(x)=o_{P}\bigg{(}\frac{1}{\sqrt{n}}\bigg{)}.\eqno(3.22)$$

 We next show that
 $$F_{n,K}(x)-\Lambda_{n,K}(x)\stackrel{d}{\longrightarrow} N(0, F(x)(1-F(x))).\eqno(3.23)$$
Note that $F_{n,K}(x)-\Lambda_{n,K}(x)$ is a martingale with respect
to an increasing $\sigma$-algebra
$\mathscr{F}_{n}=\sigma(X_{1},...,X_{n})$. So in order to verify
(3.23),  we only need to check the  conditions on the CLT for
martingales. For simplicity, set
$$E_{n,2}^{K}(x):=h_{n}^{-2}E\bigg{(}\int_{-\infty}^{x}K_{1}(t)dt\bigg{)}^{2}.\eqno(3.24)$$
We claim that the limit of $E_{n,2}^{K}(x)$ exists for any fixed $x\in\mathbb{R}$ as $n\rightarrow\infty$.  The proof is as follows: Denote
 $$G_{K}(x):=\int_{-\infty}^{x}K(u)du.$$
Recalling that $\int_{-\infty}^{\infty}K(x)dx=1$. Obviously,
$G_{K}(x)$ is the distribution function of a finite measure. Then by
the symmetry of kernel $K$, we have
$$h_{n}^{-1}\int_{-\infty}^{x}K_{1}(t)dt=h_{n}^{-1}\int_{-\infty}^{x}K\Big{(}\frac{X_{1}-t}{h_{n}}\Big{)}dt
=\int_{-\infty}^{\frac{x-X_{1}}{h_{n}}}K(u)du=G_{K}\bigg{(}\frac{x-X_{1}}{h_{n}}\bigg{)}.\eqno(3.25)$$
Note that $G_{K}$ is bounded from above by one almost surely, it then turns out that
$$E_{n,2}^{K}(x)=E\bigg{[}G_{K}\bigg{(}\frac{x-X_{1}}{h_{n}}\bigg{)}\bigg{]}^{2}\leq1.\eqno(3.26)$$
Thus we have for any fixed $x\in\mathbb{R}$,
\begin{eqnarray*}
~~~~~~~~~~~~~~~~~&&E\bigg{[}G_{K}\bigg{(}\frac{x-X_{1}}{h_{n}}\bigg{)}\bigg{]}^{2}=\int_{-\infty}^{\infty}\bigg{[}G_{K}\bigg{(}\frac{x-u}{h_{n}}\bigg{)}\bigg{]}^{2}f(u)du\\
&&~~=\bigg{(}\int_{-\infty}^{x}+\int_{x}^{\infty}\bigg{)}\bigg{\{}\bigg{[}G_{K}\bigg{(}\frac{x-u}{h_{n}}\bigg{)}\bigg{]}^{2}f(u)du\bigg{\}}
\rightarrow F(x). ~~~~~~~~~~~~~~~~~~~(3.27)
\end{eqnarray*}

The conditional variance
\begin{eqnarray*}
&&\frac{1}{nh_{n}^{2}}\sum_{i=1}^{n}E\bigg{\{}\bigg{[}\int_{-\infty}^{x}K_{i}(t)dt
 -E\bigg{(}\int_{-\infty}^{x}K_{i}(t)dt\Big{|}\mathscr{F}_{i-1}\bigg{)}\bigg{]}^{2}\bigg{|}\mathscr{F}_{i-1}\bigg{\}}\\
 &&~~=\frac{1}{nh_{n}^{2}}\sum_{i=1}^{n}E\bigg{[}\bigg{(}\int_{-\infty}^{x}K_{i}(t)dt\bigg{)}^{2}\Big{|}\mathscr{F}_{i-1}\bigg{]}
 -\frac{1}{nh_{n}^{2}}\sum_{i=1}^{n}\bigg{[}E\bigg{(}\int_{-\infty}^{x}K_{i}(t)dt\Big{|}\mathscr{F}_{i-1}\bigg{)}\bigg{]}^{2}=:\Xi^{'}-\Xi^{''}.
\end{eqnarray*}
For $\Xi^{'}$, by Lemma A.2 and Jensen$^{,}$s inequality, it follows
for any $x\in\mathbb{R}$,
\begin{eqnarray*}
~~~~~~~~~~~~~&&P\big{(}|\Xi^{'}-F(x)|>\epsilon\big{)}\\
&&~~=P\big{(}|\Xi^{'}-E_{n,2}^{K}(x)+E_{n,2}^{K}(x)-F(x)|>\epsilon\big{)}\\
&&~~\leq P\big{(}|\Xi^{'}-E_{n,2}^{K}(x)|>\epsilon/2)\leq C\epsilon^{-2}n^{-2}EG_{K}^{4}((x-X_{1})/h_{n}\big{)}\\
&&~~\leq C\epsilon^{-2}n^{-1}\rightarrow0,~~~~n\rightarrow\infty.
~~~~~~~~~~~~~~~~~~~~~~~~~~~~~~~~~~~~~~~~~~~~~~~~~~~~~~~~~(3.28)
\end{eqnarray*}
As for $\Xi^{''}$. First, similarly to that of (3.27), we have
$$E_{n,1}^{K}:=E\bigg{[}G_{K}\bigg{(}\frac{x-X_{1}}{h_{n}}\bigg{)}\bigg{]}\rightarrow F(x),~~~n\rightarrow\infty.\eqno(3.29)$$
Moreover, we have  for large $n$,
\begin{eqnarray*}
~~~~~~~~~~~~~~~~~~~~~&&P\big{(}|\Xi^{''}-F^{2}(x)|>\epsilon\big{)}\\
&&~~=P\big{(}|\Xi^{''}-(E_{n,1}^{K}(x))^{2}+(E_{n,1}^{K}(x))^{2}-F^{2}(x)|>\epsilon\big{)}\\
&&~~\leq P\big{(}|\Xi^{''}-(E_{n,2}^{K}(x))^{2}|>\epsilon/2\big{)}.
~~~~~~~~~~~~~~~~~~~~~~~~~~~~~~~~~~~~~~~~~~~~~~(3.30)
\end{eqnarray*}
Further, note that
$$\Xi^{''}-(E_{n,2}^{K}(x))^{2}=\frac{1}{n}\sum_{i=1}^{n}\bigg{[}G_{K}\bigg{(}\frac{x-X_{1}}{h_{n}}\bigg{)}-E_{n,1}^{K}(x)\bigg{]}
\bigg{[}G_{K}\bigg{(}\frac{x-X_{1}}{h_{n}}\bigg{)}+E_{n,1}^{K}(x)\bigg{]}.\eqno(3.31)$$
By the a.s. boundness of $G_{K}$, (3.29) and Lemma A.2,  the right-hand side of (3.30) is less than or equal to
$$P\bigg{(}\frac{1}{n}\sum_{i=1}^{n}\bigg{[}G_{K}\bigg{(}\frac{x-X_{1}}{h_{n}}\bigg{)}
-E_{n,1}^{K}(x)\bigg{]}>\epsilon/4\bigg{)}\rightarrow0,~~~n\rightarrow\infty.\eqno(3.32)$$
 Then it turns out that
$$\sum_{i=1}^{n}E\Big{\{}\big{[}K_{i}(x)
 -E\big{(}K_{i}(x)|\mathscr{F}_{i-1}\big{)}\big{]}^{2}|\mathscr{F}_{i-1}\Big{\}}
 \stackrel{P}{\longrightarrow}F(x)(1-F(x)).\eqno(3.33)$$
Similarly to that of (3.18),   one can show that the Lindeberg  condition holds.
Finally, by the CLT for martingales  and Slutsky$^{,}$s  theorem, we  obtains   (3.20).

  \vspace{3mm}
   Before stating  the next result, we introduce the following  condition:
  \\
  \\
   (B3) $h_{n}\searrow0$, $ nh_{n}\rightarrow\infty$  and $\sqrt{n}\omega(n)h_{n}\rightarrow0$, where $\omega(n)$
   is a nonnegative real function  satisfying $\omega(n)\nearrow\infty$ as $n\rightarrow\infty$.
  \\
  \\
  {\bf Theorem 3.3.} \textit{Suppose that condition (B3) holds, and that $\sum_{i=0}^{\infty}\rho(2^{i})<\infty$.
  Further, assume that the density function $f(x)$ is positive
   and Lipschitz continuous on $\mathbb{R}$, and that the kernel $K$ is a symmetric  function with bounded support, $\int_{-\infty}^{\infty}K(x)dx=1$
   and $\int_{-\infty}^{\infty}|K(x)|dx<\infty$.
 Then we have
 $$\sqrt{n}(F_{n,K}(x)-F(x))\stackrel{d}{\longrightarrow} N(0, F(x)(1-F(x))).\eqno(3.34)$$}
 \textit{Proof.} Observe that
 $$F_{n,K}(x)-F(x)=[F_{n,K}(x)-\Lambda_{n,K}(x)]+[\Lambda_{n,K}(x)-EF_{n,K}(x)]+[EF_{n,K}(x)-F(x)],\eqno(3.35)$$
 where
 $\Lambda_{n,K}(x)$ is defined   in the proof of Theorem 3.2. In fact, according to  Theorem 3.2, one only needs to show that
$$EF_{n,K}(x)-F(x)=o\bigg{(}\frac{1}{\sqrt{n}}\bigg{)}.\eqno(3.36)$$
 Observe that $f(x)$ is integrable on $\mathbb{R}$. Then  for any given $\epsilon>0$, it can be decomposed as follows:
 $$f(x)=f_{1}(x)+f_{2}(x),\eqno(3.37)$$
 where  $f_{1}(x)$ is continuous on a compact interval $[c_{1}, c_{2}]$ (say), and $f_{2}(x)$  satisfies
 $$\int_{-\infty}^{\infty}|f_{2}(x)|dx<\epsilon.\eqno(3.38)$$
  Denote $\tau_{n}=1/(\sqrt{n}\omega(n))$. Recalling that $\int_{-\infty}^{\infty}K(x)dx=1$, we then have for any fixed $x\in\mathbb{R}$,
 \begin{eqnarray*}
 &&|EF_{n,K}(x)-F(x)|=\bigg{|}\int_{-\infty}^{x}\bigg{[}\int_{-\infty}^{\infty}K(u)\big{[}f(t+uh_{n})-f(t)\big{]}du\bigg{]}dt\bigg{|}\\
 &&~~\leq \int_{-\infty}^{x}\int_{|u|>\tau_{n}h_{n}^{-1}}|K(u)||f_{2}(t+uh_{n})-f_{2}(t)|dudt\\
 &&~~~~~~~
 + \int_{-\infty}^{x}\int_{|u|\leq\tau_{n}h_{n}^{-1}}|K(u)||f_{2}(t+uh_{n})-f_{2}(t)|dudt=:\Theta_{1}+\Theta_{2}.
\end{eqnarray*}
Recall that the kernel $K$ is supported on the bounded interval.
Clearly,  $\Theta_{1}$ is controlled by
 \begin{eqnarray*}
~~~~~~~~~~~~~~
&&\int_{|u|>\tau_{n}h_{n}^{-1}}|K(u)|\bigg{[}\int_{-\infty}^{\infty}|f_{2}(t+uh_{n})-f_{2}(t)|dt\bigg{]}du\\
&&~~\leq2\int_{|u|>\tau_{n}h_{n}^{-1}}|K(u)|\bigg{[}\int_{-\infty}^{\infty}|f_{2}(t)|dt\bigg{]}du\rightarrow0,~~~n\rightarrow\infty.
~~~~~~~~~~~~~~~~~~~~(3.39)
\end{eqnarray*}
 For $\Theta_{2}$, note that $f_{1}$ is Lipschitz
continuous on $[c_{1}, c_{2}]$.  Subsequently,
 \begin{eqnarray*}
~~~~~~~~~~~~~~~~
 &&\Theta_{2}\leq \int_{|u|\leq\tau_{n}h_{n}^{-1}}|K(u)|
 \bigg{[}\int_{c_{1}-\tau_{n}}^{(c_{2}+\tau_{n})\wedge x}\sup_{|s|\leq\tau_{n}}|f_{1}(t+s)-f_{1}(t)|dt\bigg{]}du\\
&&~~~~\leq
 C(c_{2}-c_{1})\tau_{n}
+o\bigg{(}\frac{1}{\sqrt{n}}\bigg{)}=o\bigg{(}\frac{1}{\sqrt{n}}\bigg{)}.
~~~~~~~~~~~~~~~~~~~~~~~~~~~~~~~~~~~~~~(3.40)
\end{eqnarray*}
Therefore,  we complete the proof of (3.36).

  \sec{4. \quad Convergence rates of $\|f_{n,K}(x)-Ef_{n,K}(x)\|_{p}$  in  sup-norm   and integral $L^{p}$-norm}
   One may be interested in the consistency of   $\|f_{n,K}(x)-Ef_{n,K}(x)\|_{p}$, which are investigated  in this section.
    Among which, the uniform convergence rate with respect to  $L^{p}$-norm distance
     is established in Theorem 4.1,  while the convergence rate for integral $L^{p}$-norm   is given in Theorem 4.2.
   \\
   \\
 {\bf Theorem 4.1.} \textit{Let $p\geq2$. Suppose that $\sum_{i=0}^{\infty}\rho^{2/p}(2^{i})<\infty$,
and that  the conditions (B1), (C1), (K1) are satisfied. Then we have
$$\sup_{x\in \mathbb{R}}\big{\|}f_{n,K}(x)-Ef_{n,K}(x)\big{\|}_{p}=O\Big{[}\Big{(}\frac{1}{nh_{n}}\Big{)}^{1/2}\Big{]}.\eqno(4.1)$$}
\\
\textit{ Proof.}
 We  will use the symbols such as  $p_{k}, q_{k}, r_{k}, \mathbb{I}_{k}(m), \mathbb{J}_{k}(m)$ etc.  appeared  in the proof of Theorem 3.1. However, the values of $\alpha$ and $\beta$ are different from those as in
  the proof of Theorem 3.1. Here we select some $(p-2)/(p-1)<\alpha<1$ and $0<\beta<\min(\alpha, 1+p\alpha-p)$. In fact, $\beta$ allows taking the value $1+p\alpha-p$.  For simplicity,  define
  $$Y_{k}^{m}(x)=\sum_{i\in \mathbb{I}_{k}(m)}K_{i}(x),
  ~Z_{k}^{m}(x)=\sum_{i\in \mathbb{J}_{k}(m)}K_{i}(x),~1\leq m\leq r_{k};~Z_{k}^{r_{k}+1}(x)=\sum_{i\in \mathbb{J}_{k}(r_{k}+1)}K_{i}(x),$$
 $\xi_{k}^{m}(x)=Y_{k}^{m}(x)-E[Y_{k}^{m}(x)|\mathscr{F}_{k}(m-1)]$, where
 $$\mathscr{F}_{k}(m)=\sigma\big{(}X_{r}, r\leq 2^{k}-q_{k-1}\textrm{I}(m=1)+[(m-1)(p_{k}+q_{k})-q_{k}]\textrm{I}(m\neq1)\big{)}.$$
 Then we have for any fixed $x\in\mathbb{R}$,
 \begin{eqnarray*}
 &&f_{n,K}(x)-Ef_{n,K}(x)=\frac{1}{nh_{n}}\Bigg{\{}\Bigg{[}\sum_{k=0}^{k_{n}-1}\sum_{m=1}^{r_{k}}\xi_{k}^{m}(x)
 +\sum_{m=1}^{m_{n}-1}\xi_{k_{n}}^{m}(x)\Bigg{]}\\
  &&~~~~~~~~~~~~~~~~~~~~~~~~~~~~+ \Bigg{[}\sum_{k=0}^{k_{n}-1}\sum_{m=1}^{r_{k}+1}(Z_{k}^{m}(x)-EZ_{k}^{m}(x))+\sum_{m=1}^{m_{n}-1}(Z_{k_{n}}^{m}(x)-EZ_{k_{n}}^{m}(x))\Bigg{]}\\
  &&~~~~~~~~~~~~~~~~~~~~~~~~~~~~+ \Bigg{[}\sum_{k=0}^{k_{n}-1}\sum_{m=1}^{r_{k}}U_{k}^{m}(x)
 +\sum_{m=1}^{m_{n}-1}U_{k_{n}}^{m}(x)\Bigg{]}
+\sum_{i=N_{n}}^{n}(K_{i}(x)-EK_{i}(x))\Bigg{\}}\\
 &&~~~~~~~~~~~~~~~~~~~~~~~~~~=:\frac{1}{nh_{n}}(\textrm{I}_{1}(x)+\textrm{I}_{2}(x)+\textrm{I}_{3}(x)+\textrm{I}_{4}(x)),
  \end{eqnarray*}
where
$U_{k}^{m}(x)=E(Y_{k}^{m}(x)|\mathscr{F}_{k}(m-1))-EY_{k}^{m}(x)$.

 Thus in order to prove (4.1), it is enough to show that
 $$\sup_{x\in\mathbb{R}}\|\textrm{I}_{2}(x)+\textrm{I}_{3}(x)+\textrm{I}_{4})(x)\|_{p}=o\Big{[}(nh_{n})^{1/2}\Big{]}\eqno(4.2)$$
and
$$\sup_{x\in\mathbb{R}}\|\textrm{I}_{1}(x)\|_{p}=O\Big{[}(nh_{n})^{1/2}\Big{]}.\eqno(4.3)$$
 The proof of (4.2) will be divided into two steps:

 \textit{Step 1.} We first show for any $x\in\mathbb{R}$,
  $$ \bigg{\|}\sum_{m=1}^{r_{k}+1}(Z_{k}^{m}(x)-EZ_{k}^{m}(x))\bigg{\|}_{p}=o\Big{[}(2^{k}h_{2^{k}})^{1/2}\Big{]}.\eqno(4.4)$$
  Using Lemma A.4, and note that $\sum_{i=0}^{\infty}\rho^{2/p}(2^{i})<\infty$, it follows that
 $$E|Z_{k}^{m}(x)-EZ_{k}^{m}(x)|^{p}\leq Cq_{k}^{p/2}\|K_{1}(x)\|_{2}^{p}+Cq_{k}\|K_{1}(x)\|_{p}^{p}.\eqno(4.5)$$
 Recalling conditions (C1) and (K1), then by a simple calculation,  we have
 $$\|K_{1}(x)\|_{2}^{p}\leq (\|f\|_{\infty}\|K\|_{\infty}\|K\|_{1}h_{2^{k}})^{p/2}  \eqno(4.6)$$
 and
 $$\|K_{1}(x)\|_{p}^{p}\leq \|f\|_{\infty}\|K\|_{\infty}^{p-1}\|K\|_{1}h_{2^{k}}.  \eqno(4.7)$$
 Applying (4.5)$-$(4.7) and Minkowski$^{,}$s inequality yields
 \begin{eqnarray*}
&& \Bigg{\|}\sum_{k=0}^{k_{n}-1}\sum_{m=1}^{r_{k}+1}(Z_{k}^{m}(x)-EZ_{k}^{m}(x)) \Bigg{\|}_{p}\\
&&~~~\leq C\sum_{k=0}^{k_{n}-1}\Big{(}r_{k}h_{2^{k}}^{1/2}q_{k}^{1/2}+r_{k}h_{2^{k}}^{1/p}q_{k}^{1/p}\Big{)}\\
&&~~~=O\Big{[}(nh_{n})^{1/2}\Big{(}n^{(1+\beta-2\alpha)/2}+o(n^{1-\alpha-(1-\beta)/p})\Big{)}\Big{]}
=o\Big{[}(nh_{n})^{1/2}\Big{]},
\end{eqnarray*}
where the first equality is obtained by $1/h_{n}=o(n)$, and the
second one is due to $\beta\leq 1+p\alpha-p$. Note that
$h_{n}\searrow0$ and $nh_{n}\rightarrow\infty$ implies that
$h_{n}/h_{n+1}\leq2$ for  $n\geq1$.
 Similarly, we have
 $$\sup_{x\in\mathbb{R}}\bigg{\|}\sum_{m=1}^{m_{n}-1}(Z_{k_{n}}^{m}(x)-EZ_{k_{n}}^{m}(x))\bigg{\|}_{p}
 =o\Big{[}(nh_{n})^{1/2}\Big{]}. \eqno(4.8)$$

\textit{Step 2.} We next prove  for any $x\in\mathbb{R}$,
 $$\bigg{\|}\sum_{m=1}^{r_{k}}U_{k}^{m}(x)\bigg{\|}_{p}=o\Big{[}(2^{k}h_{2^{k}})^{1/2}\Big{]}. \eqno(4.9)$$
Recalling that $\sum_{i=0}^{\infty}\rho^{2/p}(2^{i})<\infty$, hence
w.l.o.g.,  we suppose that $\rho(n)\leq (\log n)^{-p/2}$.  By Lemma
A.4, we have
 \begin{eqnarray*}
&&E|U_{k}^{m}(x)|^{p}\leq L(\log 2p_{k})^{p}
            \bigg{[}\sum_{i=1}^{p_{k}}\rho^{2/(p-1)}(q(i/2))\bigg{]}\|K_{1}(x)\|_{p}^{p}\\
 &&~~~~~~~~~~~~~~~~~+L(\log 2p_{k})^{p}
            \bigg{[}\sum_{i=1}^{p_{k}}\rho^{2}(q(i/2))\bigg{]}^{p/2}\|K_{1}(x)\|_{2}^{p}\\
&&~~~~~~~~~~~~~~~\leq L\Big{[}k^{p}2^{-\alpha
k/(p-1)}h_{2^{k}}+k^{p}2^{\alpha(1-p)p/2}h_{2^{k}}^{p/2}\Big{]},
 \end{eqnarray*}
where $q(x)$ is the linear interpolating function of $q_{k}$.
Subsequently,
$$\|U_{k}^{m}(x)\|_{p}\leq L\Big{[}k2^{-\alpha k/p(p-1)}h_{2^{k}}^{1/p}
+k2^{\alpha(1-p)/2}h_{2^{k}}^{1/2}\Big{]}.\eqno(4.10)$$
Then taking (4.10) back into (4.9), a standard computation leads to
 \begin{eqnarray*}
&& \sup_{x\in\mathbb{R}}\bigg{\|}\sum_{k=0}^{k_{n}-1}\sum_{m=1}^{r_{k}}U_{k}^{m}(x)\bigg{\|}_{p}\\
&&~~~\leq C\sum_{k=0}^{k_{n}-1}\Big{[}r_{k}k2^{-\alpha k/p(p-1)}h_{2^{k}}^{1/p}
+r_{k}k2^{\alpha(1-p)/2}h_{2^{k}}^{1/2}\Big{]}\\
&&~~~= C\sum_{k=0}^{k_{n}-1}\Big{[}k2^{(1-\alpha-\alpha/p(p-1))k}h_{2^{k}}^{1/p}
+k2^{(1-\alpha+\alpha(1-p)/2)k}h_{2^{k}}^{1/2}\Big{]}\\
&&~~~=O\Big{[}(nh_{n})^{1/2}(\log n)\Big{(}n^{1/2-\alpha-\alpha/p(p-1)}h_{n}^{1/p-1/2}+n^{(1-\alpha-p\alpha)/2}\Big{)}\Big{]}\\
&&~~~=O\Big{[}(nh_{n})^{1/2}(\log n)
\Big{(}n^{1-\alpha-1/p}+n^{(1-\alpha-p\alpha)/2}\Big{)}\Big{]}=o\Big{[}(nh_{n})^{1/2}\Big{]},
\end{eqnarray*}
where the third equality is obtained by $nh_{n}\rightarrow\infty$,
and the last equality is due to $\alpha>(p-1)/p$.

 Similarly, we have
 $$\sup_{x\in\mathbb{R}}\bigg{\|}\sum_{m=1}^{m_{n}-1}U_{k_{n}}^{m}(x)\bigg{\|}_{p}
 =o\Big{[}(nh_{n})^{1/2}\Big{]}\eqno(4.11)$$
 and
$$\sup_{x\in\mathbb{R}}\bigg{\|}\sum_{i=N_{n}}^{n}(K_{i}(x)-EK_{i}(x))\bigg{\|}_{p}
 =o\Big{[}(nh_{n})^{1/2}\Big{]}.\eqno(4.12)$$
According to the two steps as above, we complete the proof of (4.2).

 In order to  prove  (4.3), we first show for any integer $s\geq1$,
 $$\sup_{x\in\mathbb{R}}\|\textrm{I}_{1}(x)\|_{2^{s}}=O\Big{[}(nh_{n})^{1/2}\Big{]}.\eqno(4.13)$$
 Note that
 $$\bigg{\|}\sum_{k=0}^{k_{n}-1}\sum_{m=1}^{r_{k}}\xi_{k}^{m}(x)\bigg{\|}_{2^{s}}
 \leq\sum_{k=0}^{k_{n}-1}\bigg{\|}\sum_{m=1}^{r_{k}}\xi_{k}^{m}(x)\bigg{\|}_{2^{s}}.\eqno(4.14)$$
 Hence, we only need to show that
  $$\bigg{\|}\sum_{m=1}^{r_{k}}\xi_{k}^{m}(x)\bigg{\|}_{2^{s}}=O\Big{[}(2^{k}h_{2^{k}})^{1/2}\Big{]}.\eqno(4.15)$$
  (4.15) will be derived
  by  induction on $s$: If $s=1$, using the orthogonal property of the martingale sequences and Lemma A.2, we have
$$E\bigg{[}\sum_{m=1}^{r_{k}}\xi_{k}^{m}(x)\bigg{]}^{2}=\sum_{m=1}^{r_{k}}E(\xi_{k}^{m}(x))^{2}
 \leq4\sum_{m=1}^{r_{k}}E(Y_{k}^{m}(x))^{2}=O\big{(}2^{k}h_{2^{k}}\big{)}.\eqno(4.16)$$
  Suppose that (4.15) holds true for any integer less that $s$, we
next show that it remains valid for $s$ itself. Applying  the the Marcinkiewicz$-$Zygmund$-$Burkholder inequality yields
  \begin{eqnarray*}
  &&E\bigg{[}\sum_{m=1}^{r_{k}}\xi_{k}^{m}(x)\bigg{]}^{2^{s}}\leq CE\bigg{[}\sum_{m=1}^{r_{k}}(\xi_{k}^{m}(x))^{2}\bigg{]}^{2^{s-1}}\\
  &&~~~~= CE\bigg{\{}\sum_{m=1}^{r_{k}}\Big{[}(\xi_{k}^{m}(x))^{2}-E\big{(}(\xi_{k}^{m}(x))^{2}|\mathscr{F}_{k}(m-1)\big{)}
  +E\big{(}(\xi_{k}^{m}(x))^{2}|\mathscr{F}_{k}(m-1)\big{)}\Big{]}\bigg{\}}^{2^{s-1}}\\
  &&~~~~\leq CE\bigg{\{}\sum_{m=1}^{r_{k}}\Big{[}(\xi_{k}^{m}(x))^{2}-E\big{(}(\xi_{k}^{m}(x))^{2}|\mathscr{F}_{k}(m-1)\big{)}\bigg{\}}^{2^{s-1}}\\
  &&~~~~~~~+ CE\bigg{\{}\sum_{m=1}^{r_{k}}E\big{(}(\xi_{k}^{m}(x))^{2}|\mathscr{F}_{k}(m-1)\big{)}\Big{]}\bigg{\}}^{2^{s-1}}
  =:\textrm{II}_{1}+\textrm{II}_{2},
   \end{eqnarray*}
 where the definition of $\mathscr{F}_{k}(m-1)$ can be referred  to the beginning  of the proof.

 Note that $\xi_{k}^{m}(x))^{2}-E\big{(}(\xi_{k}^{m}(x))^{2}|\mathscr{F}_{k}(m-1), m=1,2,...,$  are martingale differences.  By the induction hypothesis,  $\textrm{II}_{1}$ is of order $O\big{(}(2^{k}h_{2^{k}})^{2^{s-2}}\big{)}$.

 Using Lemma A.4 and Jensen$^{,}$s inequality,  we have
\begin{eqnarray*}
&&\textrm{II}_{2}
\leq CK(\log 2r_{k})^{2^{s-1}}
            \bigg{[}\sum_{i=1}^{r_{k}}\rho^{2^{-s}}(q(i/2))\bigg{]}\|(\xi_{k}^{m}(x))^{2}\|_{2^{s-1}}^{2^{s-1}}\\
 &&~~~~~~~+C K(\log 2r_{k})^{2^{s-1}}
            \bigg{[}\sum_{i=1}^{r_{k}}\rho^{2}(q(i/2))\bigg{]}^{2^{s-2}}\|(\xi_{k}^{m}(x))^{2}\|_{2}^{2^{s-1}}\\
&&~~~~\leq CK(\log 2r_{k})^{2^{s-1}}
            \bigg{[}\sum_{i=1}^{r_{k}}\rho^{2^{-s}}(q(i/2))\bigg{]}\|Y_{k}^{m}(x)\|_{2^{s}}^{2^{s}}\\
 &&~~~~~~~+C K(\log 2r_{k})^{2^{s-1}}
            \bigg{[}\sum_{i=1}^{r_{k}}\rho^{2}(q(i/2))\bigg{]}^{2^{s-2}}\|Y_{k}^{m}(x)\|_{4}^{2^{s}}
    =:\textrm{II}_{21}+\textrm{II}_{22}.
 \end{eqnarray*}
 Note that $\rho(n)\leq (\log n)^{-2^{s-1}}$. By Lemma A.2, a standard calculation yields
$$\textrm{II}_{21}=O\Big{[}k^{2^{s-1}}2^{(1-\alpha)k/2}\Big{(}2^{\alpha2^{s-1}k}h_{2^{k}}^{2^{s-1}}+2^{\alpha k}h_{2^{k}}\Big{)}\Big{]}=o\big{(}2^{2^{s-1}k}h_{2^{k}}^{2^{s-1}}\big{)}.\eqno(4.17)$$
 Similarly,
\begin{eqnarray*}
~~~~~~~~&&\textrm{II}_{22}=O\Big{[}k^{2^{s-1}}2^{(1-\alpha)(1-2^{s})2^{s-2}k}
\Big{(}2^{\alpha2^{s-1}k}h_{2^{k}}^{2^{s-1}}+2^{\alpha2^{s-2} k}h_{2^{k}}^{2^{s-2}}\Big{)}\Big{]}\\
&&~~~~~=o\Big{[}\big{(}2^{2^{s-1}k}h_{2^{k}}^{2^{s-1}}\big{)} \times
\big{(}2^{(\alpha-1)k}+2^{\alpha-\alpha2^{s-2}}\big{)}\Big{]}=o\big{(}2^{2^{s-1}k}h_{2^{k}}^{2^{s-1}}\big{)}.
~~~~~~~~~~~~~~~~(4.18)
\end{eqnarray*}
Combining (4.17) and (4.18) gives (4.15).

 We next  show that   (4.15) holds true   for any  $p>2$. In fact,  there exists an integer $s\geq1$ such that
 $p\in(2^{s}, 2^{s+1})$. Then it follows from Lyapunov$^{,}$s inequality,
$$\bigg{\|}\sum_{m=1}^{r_{k}}\xi_{k}^{m}(x)\bigg{\|}_{p}^{p}
\leq\bigg{\|}\sum_{m=1}^{r_{k}}\xi_{k}^{m}(x)\bigg{\|}_{2^{s}}^{2^{s+1}-p}\bigg{\|}
\sum_{m=1}^{r_{k}}\xi_{k}^{m}(x)\bigg{\|}_{2^{s+1}}^{2p-2^{s+1}}=O\Big{[}\big{(}2^{k}h_{2^{k}}\big{)}^{p/2}\Big{]}.\eqno(4.19)$$
 Finally,  we have
$$\sup_{x\in\mathbb{R}}\bigg{\|}\sum_{k=0}^{k_{n}-1}\sum_{m=1}^{r_{k}}\xi_{k}^{m}(x)\bigg{\|}_{p}=O\Big{[}(nh_{n})^{1/2}\Big{]}.\eqno(4.20)$$
Similarly, we have
$$\sup_{x\in\mathbb{R}}\bigg{\|}\sum_{m=1}^{m_{n}-1}\xi_{k_{n}}^{m}(x)\bigg{\|}_{p}=O\Big{[}(nh_{n})^{1/2}\Big{]}.\eqno(4.21)$$
According to the proof as above, we claim that (4.1) holds true.
\\
\\
{\bf Remark 3.}
   If  conditions (C1) and (K1) in Theorem 4.1 are replaced by (C3) and (K3),
   respectively,
  we have
 $$\sup_{x\in \mathbb{R}}\Big{\|}f_{n,K}(x)-f(x)\Big{\|}_{p}=O\bigg{[}\Big{(}\frac{1}{nh_{n}}\Big{)}^{1/2}+h_{n}
 \bigg{]}.\eqno(4.22)$$
\\
  \\
{\bf Theorem 4.2.} \textit{Under the conditions of Theorem 4.1,   we have for $p\geq2$,
$$\bigg{[}\int_{-\infty}^{\infty}\big{\|}f_{n,K}(x)-Ef_{n,K}(x)\big{\|}_{p}^{p}dx\bigg{]}^{1/p}=O\Big{[}\Big{(}\frac{1}{nh_{n}}\Big{)}^{1/2}\Big{]}.\eqno(4.23)$$}
\textit{Proof.} According to the decomposition as above, we only need to prove
$$\bigg{[}\int_{-\infty}^{\infty}\big{\|}\textrm{I}_{1}(x)+\textrm{I}_{2}(x)+\textrm{I}_{3}(x)+\textrm{I}_{4}(x)\big{\|}_{p}^{p}dx\bigg{]}^{1/p}=O\Big{[}\big{(}nh_{n}\big{)}^{1/2}\Big{]}.\eqno(4.24)$$
Obviously, (4.24) will be derived if one can show that
$$\sum_{j=1}^{4}\bigg{[}\int_{-\infty}^{\infty}\big{\|}\textrm{I}_{j}(x)\big{\|}_{p}^{p}dx\bigg{]}^{1/p}=O\Big{[}\big{(}nh_{n}\big{)}^{1/2}\Big{]}.\eqno(4.25)$$
In fact, the proof of (4.25) is contained  in that of Theorem 4.1.
For example, noting that
 $$\bigg{\|}\sum_{k=0}^{k_{n}-1}\sum_{m=1}^{r_{k}}\xi_{k}^{m}(x)\bigg{\|}_{p}^{p}\leq\Bigg{[}\sum_{k=0}^{k_{n}-1}\bigg{\|}\sum_{m=1}^{r_{k}}\xi_{k}^{m}(x)\bigg{\|}_{p}\Bigg{]}^{p}.\eqno(4.26)$$
Then we have
$$\bigg{[}\int_{-\infty}^{\infty}\bigg{\|}\sum_{k=0}^{k_{n}-1}\sum_{m=1}^{r_{k}}\xi_{k}^{m}(x)\bigg{\|}_{p}^{p}dx\bigg{]}^{1/p}=O\Big{[}\big{(}nh_{n}\big{)}^{1/2}\Big{]},\eqno(4.27)$$
and the order of
$\big{\|}\sum_{m=1}^{r_{k}}\xi_{k}^{m}(x)\big{\|}_{p}$  has been
obtained in Theorem 4.1. Similarly,  one can derive the order for
other terms, here they are omitted.

\renewcommand{\theequation}{5.\arabic{equation}}
 \setcounter{equation}{0}

 \sec{5. \quad Rates of strong uniform consistency  for  KDEs}
 In this section, we consider  the  a.s. convergence rates for  KDEs, which are an active subject  in probability and statistics these years.
   Among these results, Theorem 5.1 establishes the uniform rate  over a compact set,  while Theorem 5.2  gives the same
    rate over  the whole real line. It is showed that the uniform convergence rates for mixing dependent observations
   are as good as those for  i.i.d.  ones.
 \\
\\
{\bf Theorem 5.1.} \textit{Let $D$ be a compact  subset of $\mathbb{R}$. Suppose that $\rho(1)\leq1/4$ and $\sum_{i=0}^{\infty}\rho(2^{i})<\infty$, and that  the conditions (B2), (C1), (K2)  are satisfied. Then we have
$$\sup_{x\in
D}\big{|}f_{n,K}(x)-Ef_{n,K}(x)\big{|}=O_{a.s.}\Big{[}\Big{(}\frac{|\log
h_{n}|}{nh_{n}}\Big{)}^{1/2}\Big{]}.\eqno(5.1)$$}
\\
\textit{ Proof.}  We first introduce  some  notation:
    Let $\mathbb{H}_{k}$ denote the set of all integers in the interval $[2^{k}, 2^{k+1}), k\geq 0$. For $0<b<a<1/2$,
    whose values will be specified later.
     Define $p_{k}=[2^{ak}], q_{k}=[2^{bk}], r_{k}=[2^{k}/(p_{k}+q_{k})]$,  and blocks

   ~~~~~~~~~~~~~$\mathbb{I}_{k}(j)=[2^{k}+(j-1)(p_{k}+q_{k}), 2^{k}+(j-1)q_{k}+jp_{k})\cap \mathbb{N}^{+}$,

   ~~~~~~~~~~~~~$\mathbb{J}_{k}(j)=[2^{k}+(j-1)q_{k}+jp_{k}, 2^{k}+j(q_{k}+p_{k}))\cap \mathbb{N}^{+}, ~~1\leq j\leq r_{k}$,

   ~~~~~~~~~~~~~$\mathbb{J}_{k}(r_{k}+1)=[2^{k}+r_{k}(p_{k}+q_{k}), 2^{k+1})\cap \mathbb{N}^{+}.$\\
 Note that (5.1) will be derived if one can show that
  $$\max_{2^{k}\leq n<2^{k+1}}\sup_{x\in D}\big{|}f_{n,K}(x)-Ef_{n,K}(x)\big{|}=O_{a.s.}\bigg{[}\Big{(}\frac{|\log h_{2^{k}}|}{2^{k}h_{2^{k}}}\Big{)}^{1/2}\bigg{]}.\eqno(5.2)$$

  In order to prove (5.2), observe that $D$ is a compact set, so one can choose
 finite open balls with centers at $x_{1},...,x_{l_{k}}$, and radius
$d_{k}=(h_{2^{k}}^{3}|\log h_{2^{k}}|/2^{k})^{1/2}$ to cover $D$.
Obviously,  the numbers of the balls are of order
$O((2^{k}/(h_{2^{k}}^{3}|\log h_{2^{k}}|))^{1/2})$.  Further denote
the $i$th ball by $B_{i}=B(x_{i}, d_{k}), 1\leq i\leq l_{k}$,
 it is easy to see that $D\subset\bigcup_{i=1}^{l_{k}} B_{i}$.
 For simplicity, write
  $$S_{n}(x):=\frac{1}{\sqrt{h_{n}}}\sum_{i=1}^{n}\big{(}K_{i}(x)-EK_{i}(x)\big{)}.\eqno(5.3)$$
By the  Lipschitz  condition on kernel $K$, we have for any
 $x\in B_{m}$ and some $U>0$,
 $$\max_{2^{k}\leq n<2^{k+1}}\big{|}S_{n}(x)-S_{n}(x_{m})\big{|}\leq U d_{k}\max_{2^{k}\leq n< 2^{k+1}}nh_{n}^{-3/2}\leq U\sqrt{2^{k+1}|\log h_{2^{k+1}}|}.\eqno(5.4)$$
Denote $\lambda_{n}=\sqrt{n|\log h_{n}|}$,   we have for  $M>U$,
 \begin{eqnarray*}
 ~~~~~~~~~~~~~
 &&P\Big{(}\max_{2^{k}\leq n<2^{k+1}}\sup_{x\in D}|S_{n}(x)|\geq 2M\lambda_{2^{k+1}}\Big{)}\\
 &&~~\leq\sum_{m=1}^{l_{k}}\max_{1\leq m\leq l_{k}}P\Big{(}\max_{2^{k}\leq n<2^{k+1}}|S_{n}(x_{m})|\geq M\lambda_{2^{k+1}}\Big{)}~~~~~~~~~~~~~~~~~~~~~~~~~~~~~~(5.5)\\
 &&~~~~+\sum_{m=1}^{l_{k}}P\Big{(}\max_{2^{k}\leq n<2^{k+1}}\sup_{x\in B_{m}}|S_{n}(x)-S_{n}(x_{m})|\geq M\lambda_{2^{k+1}}\Big{)}.
 ~~~~~~~~~~~~~~~~~~~(5.6)
 \end{eqnarray*}
Clearly, (5.6) vanishes  on account of (5.4),  so we only need to consider (5.5).
Let us first introduce some  extra notation, define
 $$U_{k}(j, x_{m})=\frac{1}{\sqrt{h_{2^{k}}}}\sum_{i\in \mathbb{I}_{k}(j)}\big{[}K_{i}(x_{m})-EK_{i}(x_{m})\big{]},
  ~~V_{k}(j, x_{m})=\frac{1}{\sqrt{h_{2^{k}}}}\sum_{i\in \mathbb{J}_{k}(j)}\big{[}K_{i}(x_{m})-EK_{i}(x_{m})\big{]},$$
 $$1\leq j\leq r_{k};~~~V_{k}(r_{k}+1, x_{m})=\frac{1}{\sqrt{h_{2^{k}}}}\sum_{i\in \mathbb{J}_{k}(r_{k}+1)}\big{[}K_{i}(x_{m})-EK_{i}(x_{m})\big{]},
 $$
 $$\zeta_{k}(j, x_{m})=U_{k}(j, x_{m})-E\big{[}U_{k}(j, x_{m})|\mathscr{F}_{k}(j-1)\big{]},~~\mathscr{F}_{k}(j)=\sigma\big{(}X_{r}, r\leq 2^{k}+(j-1)(p_{k}+q_{k})\big{)}.$$
Then,  it is easy to give  the following  decomposition:
 \begin{eqnarray*}
 &&S_{n}(x_{m})=\sum_{j=1}^{r_{k}}\zeta_{k}(j, x_{m})+\sum_{j=1}^{r_{k}}E\big{[}U_{k}(j, x_{m})|\mathscr{F}_{k}(j-1)\big{]}\\
 &&~~~~~~~~~~~~~~+\sum_{j=1}^{r_{k}+1}E\big{[}V_{k}(j, x_{m})|\mathscr{F}_{k}(j-1)\big{]}
  =:\Sigma_{1}+\Sigma_{2}+\Sigma_{3}.
  \end{eqnarray*}

 The main ideas of the proof are as follows:
 First, note that   $\Sigma_{1}$ is a  martingale,
 in order to obtain good estimation of the  tail probability for $\Sigma_{1}$, the exponential inequality
is necessary.  Clearly,  the Freedman
 inequality is suitable for the present setting.
 However, some preliminary work is required. More precisely,  one needs  to  derive the growth rates of the bound  for the martingale differences, and the bound of the sum of the conditional variances.
Second,  we will show that $\Sigma_{2}$ and  $\Sigma_{3}$ can be negligible,   that is
 $\Sigma_{2}+\Sigma_{3}$ is of order $o_{a.s.}(\lambda_{2^{k}})$.

 The procedure as above  follows from three facts below:

 (F1)  We first show that
\textit{ $\max_{1\leq j\leq r_{k}}|\zeta_{k}(j, x_{m})|\leq p_{k}$ with probability one}. We use the following decomposition:
 \begin{eqnarray*}
&&\zeta_{k}(j, x_{m})=\frac{1}{\sqrt{h_{2^{k}}}}\sum_{i\in \mathbb{I}_{k}(j)}\big{(}A_{i}(x_{m})-E[A_{i}(x_{m})|\mathscr{F}_{k}^{i-1}]\big{)}\\
&&~~~~~~~~~~~~~~~+\frac{1}{\sqrt{h_{2^{k}}}}\sum_{i\in
\mathbb{I}_{k}(j)}\big{(}E[A_{i}(x_{m})|\mathscr{F}_{k}^{i-1}]-E[A_{i}(x_{m})|\mathscr{F}_{k}(j-1)]\big{)}=:\Sigma_{11}+\Sigma_{12},
\end{eqnarray*}
where
$$A_{i}(x_{m})=K_{i}(x_{m})-EK_{i}(x_{m}),~~~\mathscr{F}_{k}^{i}(j)=\sigma(X_{r}, r\leq 2^{k}+(j-1)(p_{k}+q_{k})+i-1).$$
Thus by Markov$^{,}$s inequality,   it follows for some $s>0$,
  $$\sum_{k=1}^{\infty}P\Big{(}\max_{1\leq j\leq r_{k}}\big{|}\zeta_{k}(j, x_{m})\big{|}> p_{k}\Big{)}\leq C\sum_{k=1}^{\infty}r_{k}p_{k}^{-s}\Big{(}E\big{|}\Sigma_{11}\big{|}^{s}+E\big{|}\Sigma_{12}\big{|}^{s}\Big{)}.\eqno(5.7)$$

  It suffices to prove the term on the right-hand side of (5.7) is finite. In fact, one only needs to estimate
   $E\big{|}\Sigma_{11}\big{|}^{s}$ and $E\big{|}\Sigma_{12}\big{|}^{s}$. For the first term, denote $B_{i,1}(x_{m})=A_{i}(x_{m})-E[A_{i}(x_{m})|\mathscr{F}_{k}^{i-1}]$. Furthermore define recursively
  $B_{i,l}(x_{m})=B_{i,l-1}^{2}(x_{m})-E[B_{i,l-1}^{2}(x_{m})|\mathscr{F}_{k}^{i-1}], l\in\mathbb{N}$.
 Clearly, $B_{i,l}(x_{m}), i=2^{k}+(j-1)(p_{k}+q_{k})+1,..., 2^{k}+(j-1)q_{k}+jp_{k}$, are martingale differences.
 We will show by the induction method that
 $$E\big{|}\Sigma_{11}\big{|}^{2^{l}}=O\Big{[}p_{k}^{2^{l-1}}\Big{]},~~~l\in\mathbb{N}^{+}.\eqno(5.8)$$
 If $l=1$, recalling the conditions on $K$ and $f$, it turns out that
$$E\Sigma_{11}^{2}\leq \frac{4}{h_{2^{k}}}\sum_{i\in \mathbb{I}_{k}(j)}EK_{i}^{2}(x_{m})=O(p_{k}).\eqno(5.9)$$
 Suppose that (5.8) holds true for $3\leq l\in\mathbb{N}^{+}$. Then using the Marcinkiewicz$-$Zygmund$-$Burkholder inequality, we have for $l\in\mathbb{N}^{+}$,
\begin{eqnarray*}
 ~~~~
 && E\big{|}\Sigma_{11}\big{|}^{2^{l}}\leq
 c_{l}h_{2^{k}}^{-2^{l-1}}E\bigg{|}\sum_{i\in
 \mathbb{I}_{k}(j)}B_{i,1}^{2}(x_{m})\bigg{|}^{2^{l-1}}\\
 &&~~\leq c_{l}h_{2^{k}}^{-2^{l-1}}E\bigg{|}\sum_{i\in
 \mathbb{I}_{k}(j)}B_{i,2}(x_{m})\bigg{|}^{2^{l-1}}+ c_{l}h_{2^{k}}^{-2^{l-1}}E\bigg{|}\sum_{i\in
 \mathbb{I}_{k}(j)}E[B_{i,1}^{2}(x_{m})|\mathscr{F}_{i-1}]\bigg{|}^{2^{l-1}}.
 ~~~~(5.10)
\end{eqnarray*}
According to the induction hypothesis, one can easily get the first
term in (5.10) is of order $O(p_{k}^{2^{l-2}})$. As for the
second term in (5.10), note that with probability one,
$$E[B_{i,1}^{2}(x_{m})|\mathscr{F}_{i-1}]\leq 4h_{2^{k}}\|K\|_{2}^{2}\|f\|_{\infty},\eqno(5.11)$$
 which is of order $O(p_{k}^{2^{l-1}})$.

 We next  show (5.8) holds true for any integer $s$
 on $(2^{l}, 2^{l+1})$. By Lyapunov$^{,}$s  inequality, it turns out for large $k$,
 \begin{eqnarray*}
~~~~~~~~~~~~~~~~~~~~~~~~~
 &&E\big{|}\Sigma_{11}\big{|}^{s}=E\big{|}\Sigma_{11}\big{|}^{(2s-2^{l+1})+(2^{l+1}-s)}\\
 &&~~\leq \big{[}E\Sigma_{11}^{2^{l}}\big{]}^{(2^{l+1}-s)/2^{l}}\big{[}E\Sigma_{11}^{2^{l+1}}\big{]}^{(s- 2^{l})/2^{l}}=O\big{[}p_{k}^{s/2}\big{]}.
 ~~~~~~~~~~~~~~~~~~~(5.12)
 \end{eqnarray*}

 As to $\Sigma_{12}$. Recalling the conditions on $K$ and $f$,  we have with probability one,
 $$E[A_{i}(x_{m})|\mathscr{F}_{k}^{i-1}]=h_{2^{k}}\int_{-\infty}^{\infty}K(u)f(x_{m}+u h_{2^{k}}|\mathscr{F}_{k}^{i-1})du=O( h_{2^{k}}).\eqno(5.13)$$
 Thus it follows  for large $k$,
 $$E\big{|}\Sigma_{12}\big{|}^{s}=O(p_{k}^{s}h_{2^{k}}^{s/2}).\eqno(5.14)$$

 On account of condition (B2), one can choose large $s$, then taking (5.12) and (5.14) back into (5.7),
 the desired result is obtained  by the Borel-Cantelli lemma.

 (F2)  We also have
$$\sum_{j=1}^{r_{k}}E\big{[}\zeta_{k}^{2}(j, x_{m})|\mathscr{F}_{k}(j-1)\big{]}
 =O_{a.s.}\big{[}2^{k}\big{]}.\eqno(5.15)$$
 Which is proved as follows:
 Observe that
$$E\big{[}\zeta_{k}^{2}(j, x_{m})|\mathscr{F}_{k}(j-1)\big{]}
=E\big{[}U_{k}^{2}(j,
x_{m})|\mathscr{F}_{k}(j-1)\big{]}-\big{(}E\big{[}U_{k}(j,
x_{m})|\mathscr{F}_{k}(j-1)\big{]}\big{)}^{2}=:\Delta_{1}-\Delta_{2}.$$

For $\Delta_{2}$,  by  the definition of $\rho-$mixing (or Lemma A.1), we have
 \begin{eqnarray*}
&&E\big{(}E\big{[}U_{k}(j, x_{m})|\mathscr{F}_{k}(j-1)\big{]}\big{)}^{2}\\
&&~~=E\big{\{}U_{k}(j, x_{m})E[U_{k}(j, x_{m})|\mathscr{F}_{k}(j-1)]\big{\}}\\
&&~~
\leq4\rho(q_{k})E\big{[}U_{k}(j, x_{m})\big{]}^{2}=O\big{[}p_{k}\rho(q_{k})\big{]}.
 \end{eqnarray*}
Using $\sum_{i=0}^{\infty}\rho(2^{i})<\infty$, it turns out that
$$\sum_{k=1}^{\infty}\sum_{j=1}^{r_{k}}E\big{(}E[U_{k}(j, x_{m})|\mathscr{F}_{k}(j-1)]\big{)}^{2}/2^{k}<\infty,\eqno(5.16)$$
which implies
$$\sum_{k=1}^{\infty}\sum_{j=r_{k-1}+1}^{r_{k}}E\big{(}E[U_{k}(j, x_{m})|\mathscr{F}_{k}(j-1)]\big{)}^{2}/2^{k}<\infty.\eqno(5.17)$$
 Applying the Borel-Cantelli lemma and the Kronecher lemma gives
$$\sum_{j=1}^{r_{n}}\big{(}E[U_{n}(j, x_{m})|\mathscr{F}_{n}(j-1)]\big{)}^{2}=o_{a.s.}(2^{n}).\eqno(5.18)$$
 Thus, the estimation for $\Delta_{2}$ is finished.

For $\Delta_{1}$, let
 $$U_{k}^{'}(j, x_{m})=U_{k}(j, x_{m})\textmd{I}(|U_{k}(j, x_{m})|\leq 2^{k/2}),~~ U_{k}^{''}(j, x_{m})=U_{k}(j, x_{m})\textmd{I}(|U_{k}(j, x_{m})|> 2^{k/2}).$$
Obviously, we have
$$E\big{[}U_{k}(j, x_{m})|\mathscr{F}_{k}(j-1)\big{]}
=E\big{[}U_{k}^{'}(j, x_{m})|\mathscr{F}_{k}(j-1)\big{]}
+E\big{[}U_{k}^{''}(j,
x_{m})|\mathscr{F}_{k}(j-1)\big{]}.\eqno(5.19)$$ Note that
\begin{eqnarray*}
~~~~~~~~~~~~
&&E\big{(} E[U_{k}^{''2}(j, x_{m})|\mathscr{F}_{k}(j-1)]\big{)}=EU_{k}^{2}(j, x_{m})\textmd{I}(|U_{k}(j, x_{m})|> 2^{k/2})\\
&&~~\leq 2^{-k\delta/2}E|U_{k}(j, x_{m})|^{2+\delta}\leq 2^{-k\delta/2}\big{[}2^{ak(2+\delta)/2}+2^{ak}h_{2^{k}}^{-\delta/2}\big{]}.
~~~~~~~~~~~~~~~~~~~(5.20)
\end{eqnarray*}
Consequently,
$$\sum_{k=1}^{\infty}\sum_{j=1}^{r_{k}}2^{-k}E\big{(} E[U_{k}^{''2}(j, x_{m})|\mathscr{F}_{k}(j-1)]\big{)}<\infty.\eqno(5.21)$$
Using the Kronecher lemma yields
$$\sum_{k=1}^{n}\sum_{j=1}^{r_{k}}E\big{(} E[U_{k}^{''2}(j, x_{m})|\mathscr{F}_{k}(j-1)]\big{)}=o_{a.s.}\big{(}2^{n}\big{)}.\eqno(5.22)$$

As for the first term on the right-hand side of (5.19),
it follows that
\begin{eqnarray*}
&&P\Big{(}\sum_{j=1}^{r_{k}}\big{[}(U_{k}^{'}(j, x_{m})-EU_{k}^{'}(j, x_{m}))|\mathscr{F}_{k}(j-1)\big{]}\geq \epsilon 2^{k}\Big{)}\\
&&~~\leq C2^{-2k}E\Big{[}\sum_{j=1}^{r_{k}}\big{[}(U_{k}^{'}(j, x_{m})-EU_{k}^{'}(j, x_{m}))|\mathscr{F}_{k}(j-1)\big{]}\Big{]}^{2}\\
&&~~\leq C2^{-2k}\sum_{j=1}^{r_{k}}\big{(}\rho(q_{k})\big{)}^{2}E[U_{k}^{'}(j, x_{m})]^{2}(\log r_{k})^{2}\\
&&~~\leq C2^{-2k}\sum_{j=1}^{r_{k}}k^{-2}2^{2ak}(\log r_{k})^{2}\leq C2^{(a-1)k}.
\end{eqnarray*}
The second inequality as above follows from Lemma A.4 with $p=2$.

Thus by the Borel-Cantelli lemma, we have
$$\sum_{j=1}^{r_{k}}\big{[}(U_{k}^{'}(j, x_{m})-EU_{k}^{'}(j, x_{m}))|\mathscr{F}_{k}(j-1)\big{]}=O_{a.s.}(2^{k}).\eqno(5.23)$$

According to (5.18), (5.22) and (5.23), the proof of (5.15) is complete. Furthermore, w.l.o.g.,  there exists some constant $V>0$, such that
$$2^{-k}\sum_{j=1}^{r_{k}}E\big{[}\zeta_{k}^{2}(j, x_{m})|\mathscr{F}_{k}(j-1)\big{]}
 \leq V ~~~~~a.s.\eqno(5.24)$$

 (F3) Since $\Sigma_{3}$ is a sum on the small blocks, we only need to
consider $\Sigma_{2}$. We have
  $$\sum_{j=1}^{r_{k}}E[U_{k}(j, x_{m})|\mathscr{F}_{k}(j-1)]=o_{a.s.}\Big{[}\sqrt{2^{k}|\log h_{2^{k}}|}\Big{]}.\eqno(5.25)$$
 Note that $2^{k}|\log h_{2^{k}}|\geq 2^{k}$ for $k\geq k_{0}$, clearly,  (5.25) follows  from  the stronger result below.
 $$\sum_{j=1}^{r_{k}}E[U_{k}(j, x_{m})|\mathscr{F}_{k}(j-1)]=o_{a.s.}\big{(}\sqrt{2^{k}}\big{)}.\eqno(5.26)$$
 Furthermore,  (5.26)  can be derived  from (5.17).

   On account of the preliminary work as above, we next consider (5.5), observe that
  $$P\Big{(}\max_{2^{k}\leq n<2^{k+1}}|S_{n}(x_{m})|\geq M\lambda_{2^{k+1}}\Big{)}
  =P\Big{(}\max_{1\leq j<2^{k}}|S_{2^{k}+j}(x_{m})|\geq M\lambda_{2^{k+1}}\Big{)}.\eqno(5.27)$$
   Then similarly  to the proof of Lemma 3.1 in Herrndorf [11],  we have
  \begin{eqnarray*}
   &&P\Big{(}\max_{1\leq j<2^{k}}|S_{2^{k}+j}(x_{m})|\geq M\lambda_{2^{k+1}}\Big{)}\\
    &&~~\leq  P\big{(}|S_{2^{k+1}}(x_{m})|\geq
    M\lambda_{2^{k+1}}/3\big{)}\\
    &&~~~~~~~\times \Big{[}1-4\rho^{2}(1)-\max_{1\leq j<2^{k}} P\big{(}|S_{2^{k+1}}(x_{m})-S_{2^{k}+j}(x_{m})|\geq M\lambda_{2^{k+1}}/3\big{)}\Big{]}^{-1}.
 \end{eqnarray*}
 By the Markov inequality, we have
 $$ \max_{1\leq j<2^{k}}P\big{(}|S_{2^{k+1}}-S_{2^{k}+j}|\geq M\lambda_{2^{k+1}}/3\big{)}\leq1/2.\eqno(5.28)$$
  Recalling  $\rho(1)\leq1/4$, which together with (5.28) yields
  $$P\Big{(}\max_{1\leq j<2^{k}}|S_{2^{k}+j}(x_{m})|\geq M\lambda_{2^{k+1}}\Big{)}\leq4 P\big{(}|S_{2^{k+1}}(x_{m})|\geq M\lambda_{2^{k+1}}/3\big{)}.\eqno(5.29)$$

 Finally, with the help of (F1)$-$(F3)  together with (5.29), then  by the  Freedman inequality (see, e.g., Lemma A.3 in the Appendix),
we have for $M\geq\sqrt{74V}$,
\begin{eqnarray*}
~~~~~~~~~~~~~~~
&&\sum_{k=1}^{\infty}\sum_{m=1}^{l_{k}}P\Big{(}\max_{2^{k}\leq n<2^{k+1}}|S_{n}(x_{m})|\geq M\lambda_{2^{k+1}}\Big{)}\\
&&~~\leq C\sum_{k=1}^{\infty}l_{k}P\Big{(}\Big{|}\sum_{j=1}^{r_{k}}\zeta_{k}(j,
x_{m})\Big{|}\geq M\lambda_{2^{k+1}}/6\Big{)}\\
&&~~\leq C\sum_{k=1}^{\infty}l_{k}\exp\bigg{\{}\frac{-M^{2}2^{k}|\log h_{2^{k}}|}{72
Mp_{k}2^{k/2}|\log h_{2^{k}}|^{1/2}+36V2^{k}}\bigg{\}}\\
&&~~\leq C\sum_{k=1}^{\infty}l_{k}h_{2^{k}}^{M^{2}/37V}<\infty.
~~~~~~~~~~~~~~~~~~~~~~~~~~~~~~~~~~~~~~~~~~~~~~~~~~~~~~~~~~~~(5.30)
\end{eqnarray*}
 Thus, applying  the Borel-Cantelli lemma yields  (5.1).
\\
\\
{\bf Remark 4.}   We  compare (5.1) with those of Peligrad [22] and
Shao [30] for mixing observations.  Peligrad [22] obtained the
following result: Let $D$ be a compact support subset of
$\mathbb{R}^{d}$ and $\{X_{n}, n\geq1\}$ be a sequence of
$\mathbb{R}^{d}$-valued $\phi-$mixing random variables with common
unknown density function $f(x)=f(x_{1},...,x_{d})$.  Suppose that
(B1) holds,  $f$ is continuous in a $\varepsilon$-neighborhood of
$D$, and   kernel $K$ satisfies:
\\
1) $K$  is a density function on $\mathbb{R}^{d}$,
\\
2) for any $x\in\mathbb{R}^{d}$, $K(x)\leq K_{1}<\infty$,
\\
3) $\|x\|^{d+1}K(x)\rightarrow0$ as $x\rightarrow\infty$,
\\
4) $\int \|x\|K(x)dx<\infty$,
\\
5) $K$ is Lipschitzian continuous of order $\gamma$ on
$\mathbb{R}^{d}$.\\
 Further, assume that
 $$\sum_{n=1}^{\infty}\phi^{1/2}(2^{n})<\infty,\eqno(5.31)$$
  then it follows that
 $$\sup_{x\in D}| f_{n}(x)-Ef_{n}(x)|=O_{a.s.}\bigg{[}\Big{(}\frac{\log^{2} n}{nh_{n}^{d}}\Big{)}^{1/2}\bigg{]}.\eqno(5.32)$$
 Especially, if~ $h_{n}=O((\log^{2} n/n)^{1/(d+2)})$, it turns out
 that
 $$\sup_{x\in D}| f_{n}(x)-f(x)|=O_{a.s.}\bigg{[}\Big{(}\frac{\log^{2} n}{n}\Big{)}^{1/(d+2)}\bigg{]}.\eqno(5.33)$$

 If the mixing rates are strengthened to $\phi(n)=O(n^{-2-d})$,
 Shao [30] obtained
 $$\sup_{x\in D}| f_{n}(x)-Ef_{n}(x)|=O_{a.s.}\bigg{[}\Big{(}\frac{\log n}{nh_{n}^{d}}\Big{)}^{1/2}\bigg{]}.\eqno(5.34)$$
  Especially, if~ $h_{n}=O((\log n/n)^{1/(d+2)})$, it turns out that
 $$\sup_{x\in D}| f_{n}(x)-f(x)|=O_{a.s.}\bigg{[}\Big{(}\frac{\log n}{n}\Big{)}^{1/(d+2)}\bigg{]}.\eqno(5.35)$$
Under conditions 3) and  4), by Bochner$-$Parzen theorem, we have
 $Ef_{n,K}(x)-f(x)=O(h_{n}).$
Note that $\phi-$mixing is contained in $\rho-$mixing,  so Theorem
5.1  holds true for $\phi-$mixing data. Obviously, the rate in (5.1)
is better than that in (5.32) under the condition (5.31); (5.34)
achieves the best possible a.s. convergence rate, however, the
mixing rate $\phi(n)=O(n^{-2-d})$ is more stronger than that in
(5.31).
\\
\\
{\bf Theorem 5.2.} \textit{ Let $X\in L^{p}$ for $p\geq2$.  Under
the conditions of Theorem 5.1,
 we have
$$\sup_{x\in \mathbb{R}}\big{|}f_{n,K}(x)-Ef_{n,K}(x)\big{|}=O_{a.s.}\Big{[}\Big{(}\frac{|\log h_{n}|}{nh_{n}}\Big{)}^{1/2}\Big{]}.\eqno(5.36)$$}
 \\
\textit{ Proof.} In order to prove (5.36), it suffices to show that
 $$\sup_{x\leq n^{3/p}}\big{|}f_{n,K}(x)-Ef_{n,K}(x)\big{|}=O_{a.s.}\Big{[}\Big{(}\frac{|\log h_{n}|}{n h_{n}}\Big{)}^{1/2}\Big{]}
 \eqno(5.37)$$
and
 $$\sup_{x>n^{3/p}}\big{|}f_{n,K}(x)-Ef_{n,K}(x)\big{|}=O_{a.s.}\Big{[}\Big{(}\frac{|\log h_{n}|}{n h_{n}}\Big{)}^{1/2}\Big{]}.
 \eqno(5.38)$$

 Along the similar proof lines as those in Theorem 5.1, one can complete the proof of (5.37). Therefore,  we only need to
consider (5.38).  Note that $K$  is bounded with compact support and
$X\in L^{p}$, using the Markov inequality, it follows that
\begin{eqnarray*}
~~~~~~~~~~~~~~~~~~~~~
&&E\Big{[}\sup_{x>n^{3/p}}\big{|}f_{n,K}(x)-Ef_{n,K}(x)\big{|}\Big{]}\leq 2h_{n}^{-1}E\Big{[}\sup_{x>n^{3/p}}\big{|}K_{i}(x)\big{|}\Big{]}\\
&&~~\leq C h_{n}^{-1}P(|X_{i}|\geq n^{3/p}/2)\leq C
n^{-3}h_{n}^{-1}, ~~~~~~~~~~~~~~~~~~~~~~~~~~~~~~~~~~~(5.39)
\end{eqnarray*}
which implies  for
some $M>0$,
$$\sum_{n=1}^{\infty}P\Big{(}\sup_{x>n^{3/p}}\big{|}f_{n,K}(x)-Ef_{n,K}(x)\big{|}\geq M(|\log h_{n}|/n h_{n})^{1/2}\Big{)}
\leq \sum_{n=1}^{\infty}\frac{1}{n^{2}}.\eqno(5.40)$$
 Thus, (5.38) is obtained by  the Borel-Cantelli lemma.

 \renewcommand{\theequation}{A.\arabic{equation}}
 \setcounter{equation}{0}

\vskip 5mm \sec{\quad Appendix}
   We list the following basic lemmas, the first one comes from Bradley and Bryc [2],  the second  can be found in  Shao [33], while the third  one is due to Freedman [7].
     \\
      \\
   {\bf Lemma A.1.}
  \textit{Let $p, ~q>1$  with $1/p+1/q=1$. Suppose that  $X\in L^{p}(\mathscr{F}_{1}^{k})$
   and $Y\in L^{q}(\mathscr{F}_{k+n}^{\infty})$ are two $\rho-$mixing random variables. Then we have
   $$|EXY-EXEY|\leq10\big{(}\rho(n)\big{)}^{\frac{2}{p}\wedge\frac{2}{q}}\|X\|_{p}\|Y\|_{q}.$$}
       \\
   {\bf Lemma A.2.}
    \textit{Let $\{X_{n}, n\geq1\}$ be a sequence of $\rho-$mixing random variables  with mean zero and
     $\|X_{i}\|_{p}<\infty$ for some $p\geq2$.
   Then  there exists a constant $L$ depending only on $p$ and $\rho(\cdot)$ such that
  for any $k\geq0, n\geq1$,
  \begin{eqnarray*}
   && E|S_{k}(n)|^{p}
  \leq Ln^{p/2}\exp \Bigg{[}L\sum_{i=0}^{[\log n]}\rho(2^{i})\Bigg{]}\max_{k<i\leq k+n}\big{\|}X_{i}\big{\|}_{2}^{p}\\
  &&~~~~~~~~~~~~~~~+Ln\exp \Bigg{[}L\sum_{j=0}^{[\log n]}\rho^{2/p}(2^{j})\Bigg{]}\max_{k<i\leq k+n}\big{\|}X_{i}\big{\|}_{p}^{p}.
  \end{eqnarray*}}
 \\
 \\
  {\bf Lemma A.3.}
  \textit{Let $\{X_{n}, \mathscr{F}_{n},  n\geq1\}$ be a martingale difference sequence with  $S_{n}=\sum_{i=1}^{n}X_{i}$. Suppose that $\tau$ is a stopping time,
  and $L$ a positive real number. Suppose $P(|X_{i}|\leq L, i\leq \tau)=1$.
  Then for all positive real numbers $a$ and $b$,
 $$P(S_{n}\geq a,  T_{n}\leq b ~~for~ some~~ n\leq\tau)\leq\exp\bigg{[}\frac{-a^{2}}{2(La+b)}\bigg{]},$$
  where  $T_{n}=\sum_{i=1}^{n}\textrm{Var}(X_{i}|\mathscr{F}_{i-1})$.}

 \vskip 3mm
 Before formulating  the next Lemma, we  give some extra notation: Let $p_{k}$ and $q_{k}$ be  sequences of positive integers satisfying
   $q_{k}=o(p_{k})$ and $q_{k}\nearrow\infty$ as
  $k\rightarrow\infty$.   Then the successively blocks which only include integers are defined as follows:

  ~~~~~~~~~~~~~$\mathbb{I}_{k}=[(k-1)(p_{k}+q_{k})+1, (k-1)q_{k}+kp_{k})\cap \mathbb{N}^{+}$,

   ~~~~~~~~~~~~~$\mathbb{J}_{k}=[(k-1)q_{k}+kp_{k}, k(p_{k}+q_{k}))\cap
   \mathbb{N}^{+},~~k=1,2,...$\\
 Furthermore, let $\{X_{n}, n\geq1\}$ be a  sequence of $\rho-$mixing random variables with mean zero, and $E|X_{n}|^{p}<\infty$
  for some $p\geq2$.
 Denote
$$\xi_{k}=\sum_{i\in\mathbb{I}_{k}}X_{i},~~~~\eta_{k}=\sum_{i\in\mathbb{J}_{k}}X_{i}.$$
  \\
  \\
  {\bf Lemma A.4.}
 \textit{Let $\{\xi_{n}, n\geq1\}$ be as above.  Suppose that $\sum_{i=0}^{\infty}\rho(2^{i})<\infty$, then for $p\geq2$, there exists a positive constant $L=L(p, \rho(\cdot))$ depending only on $p$ and $\rho(\cdot)$ such that for every $k\geq0, n\geq1$,
  \begin{eqnarray}
 && E|G_{m}(n)|^{p}\leq L(\log 2n)^{p}
            \Bigg{[}\sum_{k=m+1}^{m+n}\rho^{2}(q(k/2))\big{\|}\xi_{k}\big{\|}_{2}^{2}\Bigg{]}^{p/2}\nonumber \\
  &&~~~~~~~~~~~~~~~~+L(\log 2n)^{p}
\Bigg{[}\sum_{k=m+1}^{m+n}\rho^{2/(p-1)}(q(k/2))\big{\|}\xi_{k}\big{\|}_{p}^{p}
\Bigg{]},
  \end{eqnarray}
 where $G_{m}(n)=\sum_{k=m+1}^{m+n}E(\xi_{k}|\mathscr{F}_{k-1})$,  $\mathscr{F}_{k}=\sigma( X_{i}, i\leq (k-1)(p_{k}+q_{k})+p_{k})$,
  and $q(x)$ is the linear interpolating function
of $q_{k}$.}\\
 \textit{Proof.} We claim that $p$ may  not be integer. However,
    we only consider the integral case in the present article
  since the proof  for the non-integral  part is similar.   We first show that (A.1) holds true for all even numbers, then verify that it is the case
  for all odd numbers.  The proof is decomposed
  into the following three  steps  by induction on $p$.
\\
\textit{Step 1.}  If $p=2$, (A.1) can be rewrote as follows:
  \begin{equation}
  EG_{m}^{2}(n)\leq L(\log 2n)^{2}
            \Bigg{[}\sum_{k=m+1}^{m+n}\rho^{2}(q(k/2))\big{\|}\xi_{k}\big{\|}_{2}^{2}\Bigg{]}.
  \end{equation}
A similar result as (A.2) can be found in Shao [31].  However, for
the reader$^{,}$s convenience, we give its proof by induction on $n$
below:   If $n=1$, by Lemma A.1, it follows that
 \begin{eqnarray}
 &&EG_{m}^{2}(1)=E(E(\xi_{m+1}|\mathscr{F}_{m}))^{2}\nonumber\\
 &&~~~~~~~~~~~=E(\xi_{m+1}E(\xi_{m+1}|\mathscr{F}_{m}))\nonumber\\
 &&~~~~~~~~~~~\leq10\rho(q_{m})\|\xi_{m+1}\|_{2}\|E(\xi_{m+1}|\mathscr{F}_{m})\|_{2}.
 \end{eqnarray}
Then a simple calculation leads to
\begin{eqnarray}
 EG_{m}^{2}(1)\leq100\rho^{2}(q_{m})\|\xi_{m+1}\|_{2}^{2}.
\end{eqnarray}

 Suppose that (A.2) holds true
for any integer less than $n$.  We next show it remains valid for
$n$ itself. Let $n_{1}=[n/2], n_{2}=n-n_{1}$. Clearly,
 \begin{eqnarray}
 &&EG_{m}^{2}(n)=EG_{m}^{2}(n_{1})+EG_{m+n_{1}}^{2}(n_{2})+2EG_{m}(n_{1})G_{m+n_{1}}(n_{2})\nonumber\\
 &&~~~~~~~~~~~\leq EG_{m}^{2}(n_{1})+EG_{m+n_{1}}^{2}(n_{2})+2\rho(q_{m+n_{1}})\|G_{m}(n_{1})\|_{2}\|G_{m+n_{1}}(n_{2})\|_{2}.
 \end{eqnarray}
 By the induction hypothesis and Lemma A.2 with $\sum_{i=0}^{\infty}\rho(2^{i})<\infty$, we
 have
 \begin{eqnarray}
&&EG_{m}^{2}(n)\leq L(\log 2n_{1})^{2}
            \Bigg{[}\sum_{k=m+1}^{m+n_{1}}\rho^{2}(q(k/2))\big{\|}\xi_{k}\big{\|}_{2}^{2}\Bigg{]}\nonumber\\
 &&~~~~~~~~~~~~~~+ L(\log 2n_{2})^{2}
            \Bigg{[}\sum_{k=m+n_{1}+1}^{m+n}\rho^{2}(q(k/2))\big{\|}\xi_{k}\big{\|}_{2}^{2}\Bigg{]}\nonumber\\
            &&~~~~~~~~~~~~~~+L\log (2n_{1})\rho(q_{m+n_{1}})\Bigg{[}\sum_{k=m+n_{1}+1}^{m+n}\big{\|}\xi_{k}\big{\|}_{2}^{2} \Bigg{]}^{1/2}
 \Bigg{[}\sum_{k=m+1}^{m+n_{1}}\rho^{2}(q(k/2))\big{\|}\xi_{k}\big{\|}_{2}^{2} \Bigg{]}^{1/2}\nonumber\\
 &&~~~~~~~~~~~\leq L(\log 2n)^{2}
   \Bigg{[}\sum_{k=m+1}^{m+n}\rho^{2}(q(k/2))\big{\|}\xi_{k}\big{\|}_{2}^{2}\Bigg{]}.
 \end{eqnarray}
\textit{Step 2.}  Let $p=2l$ for $l\geq2$.  Suppose that (A.1) holds
true for all even numbers less than $p$.  We will show that it is
also valid for $p$. To this end, the following preliminary work is
needed.

 (i) We will derive an upper bound for  the $p$th moment of
$G_{k}(n)$. The  basic inequalities  below are useful: For any
$x\geq0$ and $p>1$, we have
 \begin{eqnarray}
 (1+x)^{p}\leq 1+x^{p}+4^{p}(x+x^{p-1}).
 \end{eqnarray}
Moreover, let $\alpha,~\beta>0$ and $\alpha+\beta=1$. Applying
Young$^{,}$s inequality, we have for any $x,~y>0$,
 \begin{eqnarray}
x^{\alpha}y^{\beta}\leq \alpha x+\beta y\leq x+y.
 \end{eqnarray}
 Thus  by (A.7), it follows for $n\geq2$,
 \begin{eqnarray}
 &&EG_{m}^{p}(n)= E(G_{m}(n_{1})+G_{m+n_{1}}(n_{2}))^{p}\nonumber\\
 &&~~~~~~~~~~~\leq E(G_{m}(n_{1}))^{p}+E(G_{m+n_{1}}(n_{2}))^{p}\nonumber\\
 &&~~~~~~~~~~~~~~+4^{p}\big{(}EG_{m}(n_{1})(G_{m+n_{1}}(n_{2}))^{p-1}
 +E(G_{m}(n_{1}))^{p-1}G_{m+n_{1}}(n_{2})\big{)}.
 \end{eqnarray}
 Using Lemma A.1 and (A.8),   we give the following stronger  upper
 bound,
 \begin{eqnarray*}
 &&EG_{m}(n_{1})(G_{m+n_{1}}(n_{2}))^{p-1}\\
&&~\leq E|G_{m}(n_{1})|E|G_{m+n_{1}}(n_{2})|^{p-1}+10\rho^{2/p}(q_{m+n_{1}})\|G_{m}(n_{1})\|_{p}\|G_{m+n_{1}}(n_{2})\|_{p}^{p-1}\\
&&~\leq
E|G_{m}(n_{1})|E|G_{m+n_{1}}(n_{2})|^{p-1}+10\rho^{2/p}(q_{m+n_{1}})\Big{[}\|G_{m}(n_{1})\|_{p}^{p}+\|G_{m+n_{1}}(n_{2})\|_{p}^{p}\Big{]}.
 \end{eqnarray*}
Similarly, we have
 \begin{align*}
&E(G_{m}(n_{1}))^{p-1}G_{m+n_{1}}(n_{2})\\
&~~~\leq E|G_{m}(n_{1})|^{p-1}E|G_{m+n_{1}}(n_{2})|+
10\rho^{2/p}(q_{m+n_{1}})\Big{[}\|G_{m}(n_{1})\|_{p}^{p}+\|G_{m+n_{1}}(n_{2})\|_{p}^{p}\Big{]}.
 \end{align*}
Hence combining the estimations as above, we get
 \begin{eqnarray*}
 && EG_{m}^{p}(n)\leq\big{(}1+20\times4^{p}\rho^{2/p}(q_{m+n_{1}})\big{)}\Big{[}EG_{m}^{p}(n_{1})+EG_{m+n_{1}}^{p}(n_{2})\Big{]}\\
 &&~~~~~~~~~~~~~~+E|G_{m}(n_{1})|^{p-1}E|G_{m+n_{1}}(n_{2})|+E|G_{m}(n_{1})|E|G_{m+n_{1}}(n_{2})|^{p-1}\\
&&~~~~~~~~~~~\leq\big{(}1+20\times4^{p}\rho^{2/p}(q_{m+n_{1}})\big{)}\Big{[}E|G_{m}(n_{1})|^{p}+E|G_{m+n_{1}}(n_{2})|^{p}\Big{]}\\
  &&~~~~~~~~~~~~~~~+\parallel G_{m}(n_{1})\parallel_{p-1}^{p-1}\parallel G_{m+n_{1}}(n_{2})\parallel_{2}
  +\parallel G_{m}(n_{1})\parallel_{2}\parallel
  G_{m+n_{1}}(n_{2})\parallel_{p-1}^{p-1}.
 \end{eqnarray*}

 (ii)  We next show for each $n\geq1$,
  \begin{eqnarray}
  &&EG_{m}^{p}(n)\leq L(\log 2n)^{p}
            \Bigg{[}\sum_{k=m+1}^{m+n}\rho^{2}(q(k/2))\|\xi_{k}\|_{2}^{2}\Bigg{]}^{p/2}\nonumber\\
  &&~~~~~~~~~~~~~~+L(\log 2n)^{p}\Bigg{[}\sum_{k=m+1}^{m+n}\rho^{2/(p-2)}(q(k/2))\|\xi_{k}\|_{p-1}^{p-1}\Bigg{]}
 \Bigg{[}\sum_{k=m+1}^{m+n}\rho^{2}(q(k/2))\|\xi_{k}\|_{2}^{2}\Bigg{]}^{1/2}\nonumber\\
 &&~~~~~~~~~~~~~~+L(\log
 2n)^{p}\Bigg{[}\sum_{k=m+1}^{m+n}\rho^{2/(p-1)}(q(k/2))\|\xi_{k}\|^{p}_{p}\Bigg{]}.
  \end{eqnarray}

To prove (A.10), we apply the  induction method on $n$: If $n=1$, it
 follows from Lemma A.1,
 \begin{eqnarray}
EG_{m}^{p}(1)=E\big{[}G_{m}(1)(G_{m}(1))^{p-1}\big{]}\leq10\rho^{2/p}(q_{m})\parallel
\xi_{m+1}\parallel_{p}\parallel G_{m}(1)\parallel_{p}^{p-1}.
 \end{eqnarray}
Then a standard calculation leads to
 \begin{eqnarray}
 EG_{m}^{p}(1)\leq10\rho^{2}(q_{m})\parallel\xi_{m+1}\parallel_{p}^{p}.
 \end{eqnarray}
 Observe that $\rho^{2}(q_{m})\leq\rho^{2/(p-1)}(q_{m})$ for $l\geq1$, it is easy
 to see that (A.10) holds true.

 Suppose that (A.10) is valid for any integer less than $n$, we next
show that it is the case for $n$ itself. On account of the induction
hypothesis, $EG_{m}(n)^{p}$ is less than or equal to
  \begin{eqnarray*}
  &&L(\log 2n)^{p}
             \Bigg{\{}\Bigg{[}\sum_{k=m+1}^{m+n_{1}}\rho^{2}(q(k/2))\|\xi_{k}\|_{2}^{2}\Bigg{]}^{p/2}
              +\Bigg{[}\sum_{k=m+n_{1}+1}^{m+n}\rho^{2}(q(k/2))\|\xi_{k}\|_{2}^{2}\Bigg{]}^{p/2} \Bigg{\}}\\
&&~~~~~+L(\log 2n)^{p}
\Bigg{\{}\Bigg{[}\sum_{k=m+1}^{m+n_{1}}\rho^{2/(p-2)}(q(k/2))\|\xi_{k}\|_{p-1}^{p-1}\Bigg{]}\Bigg{[}\sum_{k=m+1}^{m+n_{1}}\rho^{2}(q(k/2))\|\xi_{k}\|_{2}^{2}\Bigg{]}^{1/2}\\
&&~~~~~~~~~~~~~~~~~~~~~~+\Bigg{[}\sum_{k=m+n_{1}+1}^{m+n}\rho^{2/(p-2)}(q(k/2))\|\xi_{k}\|_{p-1}^{p-1}\Bigg{]}\Bigg{[}\sum_{k=m+n_{1}+1}^{m+n}\rho^{2}(q(k/2))\|\xi_{k}\|_{2}^{2}\Bigg{]}^{1/2}\Bigg{\}}\\
&&~~~~~+\Bigg{[}L\sum_{k=m+1}^{m+n_{1}}\rho^{2/(p-1)}(q(k/2))\|\xi_{k}\|_{p}^{p}
+L\sum_{k=m+n_{1}+1}^{m+n}\rho^{2/(p-1)}(q(k/2))\|\xi_{k}\|_{p}^{p}\Bigg{]}\\
&&~~~~~+4^{3p}L\Bigg{\{}
\log^{1/2} (2n)\Bigg{(}\sum_{k=k+1}^{k+n}\rho^{2}(q(k/2))\|\xi_{k}\|_{2}^{2}\Bigg{)}^{1/2}(\log 2n)^{p-1}\\
&&~~~~~~~~~~~~~~~~~~\times
\Bigg{\{}\Bigg{[}\sum_{k=m+1}^{m+n}\rho^{2/(p-2)}(q(k/2))\|\xi_{k}\|_{p-1}^{p-1}\Bigg{]}
 +\Bigg{[}\sum_{k=m+1}^{m+n}\rho^{2}(q(k/2))\|\xi_{k}\|_{2}^{2}\Bigg{]}^{\frac{p-1}{2}}\Bigg{\}}\Bigg{\}}\\
  &&~~\leq L(\log 2n)^{p}
            \Bigg{[}\sum_{k=m+1}^{m+n}\rho^{2}(q(k/2))\|\xi_{k}\|_{2}^{2}\Bigg{]}^{p/2}\\
&&~~~~~+L(\log
2n)^{p}\Bigg{[}\sum_{k=m+1}^{m+n}\rho^{2/(p-1)}(q(k/2)))\|\xi_{k}\|_{p}^{p}\Bigg{]}
\\
 &&~~~~~+L(\log 2n)^{p}
\Bigg{[}\sum_{k=m+1}^{m+n}\rho^{2/(p-2)}(q(k/2))\|\xi_{k}\|_{p-1}^{p-1}\Bigg{]}\Bigg{[}\sum_{k=m+1}^{m+n}\rho^{2}(q(k/2))\|\xi_{k}\|_{2}^{2}\Bigg{]}^{1/2}.
\end{eqnarray*}
Therefore the proof of (A.10) is complete.

  (iii) We finally verify that  (A.1) holds true for $p=2l$.
 Observe that
 \begin{eqnarray}
&&(\log
2n)^{p}\Bigg{[}\sum_{k=m+1}^{m+n}\rho^{2/(p-2)}(q(k/2))\|\xi_{k}\|_{p-1}^{p-1}\Bigg{]}
\Bigg{[}\sum_{k=m+1}^{m+n}\rho^{2}(q(k/2))\|\xi_{k}\|_{2}^{2}\Bigg{]}^{1/2}\nonumber\\
&&~~\leq(\log
2n)^{p}\Bigg{[}\sum_{k=m+1}^{m+n}\rho^{2/(p-1)}(q(k/2))\|\xi_{k}\|_{p-1}^{p-1}\Bigg{]}
\Bigg{[}\sum_{k=m+1}^{m+n}\rho^{2}(q(k/2))\|\xi_{k}\|_{2}^{2}\Bigg{]}^{1/2}.
 \end{eqnarray}
Then by  Lyapunov$^{,}$s inequality, it follows that
 \begin{eqnarray}
\sum_{k=m+1}^{m+n}\rho^{2/(p-1)}(q(k/2))\|\xi_{k}\|_{p-1}^{p-1} \leq
\sum_{k=m+1}^{m+n}\rho^{2/(p-1)}(q(k/2))\|\xi_{k}\|_{2}^{2/(p-2)}\|\xi_{k}\|_{p}^{p(p-3)/(p-2)}.
 \end{eqnarray}
 Furthermore  note that
\begin{eqnarray}
\frac{2}{p-1}=\frac{2}{p-2}+\frac{2(p-3)}{(p-1)(p-2)}.
 \end{eqnarray}
Therefore  by the H\"{o}lder inequality,  the right hand side of
(A.14) is less than or equal to
\begin{eqnarray}
\Bigg{[}\sum_{k=m+1}^{m+n}\rho^{2}(q(k/2))\|\xi_{k}\|_{2}^{2}\Bigg{]}^{1/(p-2)}\Bigg{[}\sum_{k=m+1}^{m+n}\rho^{2/(p-1)}(q(k/2))\|\xi_{k}\|_{p}^{p}\Bigg{]}^{(p-3)/(p-2)}.
 \end{eqnarray}
 On account of (A.14) and  (A.16), the right hand side of (A.13) is
 controlled by
\begin{eqnarray}
&&(\log
2n)^{p}\Bigg{[}\sum_{k=m+1}^{m+n}\rho^{2}(q(k/2))\|\xi_{k}\|_{2}^{2}\Bigg{]}^{\frac{1}{2}+\frac{1}{p-2}}\Bigg{[}\sum_{k=m+1}^{m+n}\rho^{2/(p-1)}(q(k/2))\|\xi_{k}\|_{p}^{p}\Bigg{]}^{\frac{p-3}{p-2}}\nonumber\\
&&~~=(\log
2n)^{p}\Bigg{[}\sum_{k=m+1}^{m+n}\rho^{2}(q(k/2))\|\xi_{k}\|_{2}^{2}\Bigg{]}^{\frac{p}{2}\times\frac{1}{p-2}}\Bigg{[}\sum_{k=m+1}^{m+n}\rho^{2/(p-1)}(q(k/2))\|\xi_{k}\|_{p}^{p}\Bigg{]}^{\frac{p-3}{p-2}}\nonumber
\\
&&~~\leq (\log
2n)^{p}\Bigg{\{}\Bigg{[}\sum_{k=m+1}^{m+n}\rho^{2}(q(k/2))\|\xi_{k}\|_{2}^{2}\Bigg{]}^{p/2}+\sum_{k=m+1}^{m+n}\rho^{2/(p-1)}(q(k/2))\|\xi_{k}\|_{p}^{p}\Bigg{\}}.
 \end{eqnarray}

 \textit{Step 3.}   Assume that $p$ is an odd number. Clearly,  there exists integer $l\geq2$ such that $2l-2<p<2l$.
 Again using Lyapunov$^{,}$s inequality and
(A.8), we have
 \begin{eqnarray}
  &&E|G_{m}(n)|^{p}=E|G_{m}(n)|^{2l-1}\\
   &&~~~~~~~~~~~~~~\leq \big{(}EG_{m}(n)^{2l-2}\big{)}^{1/2} \big{(}EG_{m}(n)^{2l}\big{)}^{1/2}\nonumber\\
  &&~~~~~~~~~~~~~~\leq EG_{m}(n)^{2l-2}+EG_{m}(n)^{2l}.
 \end{eqnarray}

According to the procedures as above,  the proof of (A.1) is
complete.

\vskip 6mm

\vskip0.2in \no {\bf References} \vskip0.1in

 \footnotesize

\REF{[1]} Bickel, P. J., and Rosenblatt, M. (1973).  On some global
measures of the deviations of density function estimates. Ann.
Statist. 1, 1071-1095.

\REF{[2]} Bradley, R. C., and Bryc, W. (1985).  Multilinear forms and measures of dependence between random variables. J. Multivariate Anal. 16, 335-367.
measures of the deviations of density function estimates. Ann.
Statist. 1, 1071-1095.

\REF{[3]} Deheuvels, P. (2000). Uniform limit laws for kernel
density estimators on possibly unbounded intervals. In Recent
Advances in Reliability Theory: Methodology, Practice and Inference
(N. Limnios and M. Nikulin, eds.) 477-492, Birkh\"{a}user, Boston.

\REF{[4]}Einmahl, U., and Mason, D. M. (2000).  An empirical process
approach to the uniform consistency of kernel-type function
estimators.  J. Theoret. Probab. 13, 1-37.

\REF{[5]} Einmahl, U., and Mason, D. M. (2005). Uniform in bandwidth
consistency of kernel-type function estimators. Ann. Statist. 33,
1380-1403.

\REF{[6]} F\"{o}ldes, A. (1974).  Density estimation for dependent
sample.  Studia Scientiarum Mathematicarum Hungarica  9, 443-452.

\REF{[7]} Freedman, D. A. (1975). On tail probabilities for martingales.  Ann. Probab. 3, 100-118.

\REF{[8]} Gin\'{e}, E., and  Guillou, A. (2002). On consistency of
kernel density estimators for randomly censored data: Rates holding
uniformly over adaptive intervals.   Ann. Inst. H. Poincar\'{e}
Probab. Statist. 37, 503-522.

\REF{[9]} Gin\'{e}, E., and Guillou, A. (2002). Rates of strong
consistency for multivariate kernel density estimators. Ann. Inst.
H. Poincar\'{e} Probab. Statist,  38, 907-922.

\REF{[10]} Gin\'{e}, E., Koltchinskii, V., and Zinn, J. (2004).
Weighted uniform consistency of kernel density estimators.  Ann.
Probab. 32, 2570-2605.

\REF{[11]} Herrndorf, N. (1983). The Invariance Principle for p-mixing Sequences.
Z. Wahr.  verw. Gebiete  63, 97-108.

 \REF{[12]} Kolmogorov, A. N., and  Rozanov, U. A. (1960). On the strong mixing conditions of a stationary Gaussian process.  Probab Theory  Appl.  2, 222-227.

\REF{[13]} Liebscher, E. (1995). Strong convergence of sums of
$\varphi-$mixing random variables.  Math. Methods Statist. 4, 216-229.

\REF{[14]} Liebscher, E.  (1996). Strong convergence of sums of
$\alpha-$mixing random variables with applications to density
estimation.   Stochastic Process. Appl.  65,  69-80.

\REF{[15]} Lin, Z., and Lu, C. (1997).  Limit Theory on Mixing Dependent Random Variables, Science Press and Kluwer Academic Publishers.

 \REF{[16]} Neumann, M. H. (1998). Strong approximation of density estimators from weakly dependent observations by density
estimators from independent observations. Ann. Statist. 26, 2014-2048.

\REF{[17]} Nze, P. A., and Rios, R. (1995). Density estimation in the
$L^{\infty}$ norm for mixing processes.  C. R. Acad. Sci.
Paris 320, 1259-1262.

\REF{[18]} Parzen, E. (1962). On the estimation of a probability
density function and the mode.
 Ann. Math. Statist. 33, 1065-1076.

 \REF{[19]} Peligrad, M. (1982). Invariance principles for mixing sequences of random variables.
 Ann. Probab. 10, 968-981.

 \REF{[20]} Peligrad, M. (1986). Recent advances in the central limit theorem and its weak invariance principle for
 mixing sequences of random variables.  Dependence in Probab. Statist. Eberlin, E. and Taqqu, M. S. (eds)
 Progress in Probab. Statist. Birckhauser, 11, 193-223.

\REF{[21]} Peligrad, M. (1987). On the central limit theorem for $\rho-$mixing sequences of random variables.
 Ann. Probab. 15, 1387-1394.

\REF{[22]}Peligrad, M. (1992). Properties of uniform consistency of
the kernel estimators of density and of regression functions under
dependent assumptions.  Stoch. Stoch. Reports  40, 147-168.

 \REF{[23]}Peligrad, M., and Shao, Q. M. (1995). Estimation of the variance of partial sums for  $\rho-$mixing  random variables.
 J. Multivariate Anal. 52, 140-157.

 \REF{[24]} Peligrad, M., and Utev,  S. (2005). A new maximal inequality and invariance principle for stationary sequences.
          Ann.  Probab.   33, 798-815.

\REF{[25]} Rosenblatt,  M. (1956). Remarks on some nonparametric
estimates of a density function.  Ann. Math. Statist. 27, 832-835.

\REF{[26]} Roussas, G. G.  (1988). nonparametric  estimation in
mixing sequences of random variables.
 J.  Statist. Plann. Inference 18, 135-149.

\REF{[27]} R\"{u}schendorf,  L. (1977). Consistency of estimators for
multivariate density functions and for the mode.  Sankhya 39, 243-250.

\REF{[28]}Sarda, P., and  Vieu,  P. (1989).  Empirical
distribution function for mixing random variables.  Statistics 20,
559-571.

\REF{[29]} Shao, Q. M. (1989). On the invariance principle for
$\rho-$mixing sequences of random variables.  Chin. Ann. Math. 10B,
427-433.

 \REF{[30]} Shao, Q. M. (1990). Exponential inequalities and density estimation under dependent assumptions. Manuscript.

 \REF{[31]} Shao, Q. M. (1993). On the invariance principle for stationary $\rho-$mixing sequence with infinite variance. Chin. Ann. Math.
14B, 27-42.

 \REF{[32]}Shao, Q. M. (1993).  Almost sure invariance principles for mixing sequences of random variables. Stochastic Process Appl. 48, 319-334.

\REF{[33]} Shao, Q. M.  (1995).  Maximal inequality for partial sums of $\rho-$mixing sequences.  Ann. Probab. 23, 948-965.

\REF{[34]} Silverman, B. W. (1978).  Weak and strong uniform
consistency of the kernel estimate of a density and its derivatives.
Ann. Statist. 6, 177-184.

\REF{[35]} Stute,  W. (1982). A law of the logarithm for kernel
density estimators,   Ann. Probab. 10, 414-422.

\REF{[36]} Stute, W. (1984). The oscillation behavior of empirical
processes: The multivariate case.   Ann. Probab. 12, 361-379.

\REF{[37]} Woodroofe,  M. (1967). On the maximum deviation of the
sample density.  Ann. Math. Statist. 38, 475-481.

\REF{[38]}Woodroofe, M. (1970). Discussion of "Density estimates
and Markov sequences" by M. Rosenblatt. In nonparametric Techniques
in Statistical Inference (M. Puri, ed.)  211-213,  Cambridge Univ.
Press.

\REF{[39]} Wu, W. B., Huang, Y., and Huang, Y. (2010). Kernel estimation for time series: An asymptotic theory.
 Stochastic Process. Appl. 120, 2412-2431.

\REF{[40]} S. Yakowitz,   nonparametric density estimation,
prediction and regression for Markov sequences,  J. Amer. Statist.
Assoc.  80 (1985)  215-221.

\REF{[41]}B. Yu,   Density estimation in the $L^{\infty}$
norm for dependent data with applications to the Gibbs sampler,  Ann.
Statist. 21 (1993)  711-735.

\end{document}